\documentclass{article}
\usepackage{graphicx}
\usepackage{subcaption}
\usepackage{amsmath}
\usepackage{float}

\begin{document}

\title{Anomalous Diffusion:\\
Fractional Brownian Motion vs. \\
Fractional Ito Motion}
\author{Iddo Eliazar\thanks{%
E-mail: \emph{eliazar@tauex.tau.ac.il}} \and Tal Kachman\thanks{Donders Centre for Cognition, Radboud university, Netherlands.}}
\maketitle

\begin{abstract}
Generalizing Brownian motion (BM), fractional Brownian motion (FBM) is a
paradigmatic selfsimilar model for anomalous diffusion. Specifically,
varying its Hurst exponent, FBM spans: sub-diffusion, regular diffusion, and
super-diffusion. As BM, also FBM is a symmetric and Gaussian process, with a
continuous trajectory, and with a stationary velocity. In contrast to BM,
FBM is neither a Markov process nor a martingale, and its velocity is
correlated. Based on a recent study of selfsimilar Ito diffusions, we explore an alternative selfsimilar model for anomalous diffusion: \emph{%
fractional Ito motion} (FIM). The FIM model exhibits the same Hurst-exponent
behavior as FBM, and it is also a symmetric process with a continuous
trajectory. In sharp contrast to FBM, we show that FIM: is not a Gaussian
process; is a Markov process; is a martingale; and its velocity is not
stationary and is not correlated. On the one hand, FBM is hard to simulate,
its analytic tractability is limited, and it generates only a Gaussian
dissipation pattern. On the other hand, FIM is easy to simulate, it is
analytically tractable, and it generates non-Gaussian
dissipation patterns. Moreover, we show that FIM has an intimate linkage to
diffusion in a logarithmic potential. With its compelling properties, FIM
offers researchers and practitioners a highly workable analytic model for
anomalous diffusion.

\ \ 

\ 

\ 


\bigskip\ 

Brownian motion)
\end{abstract}

\newpage

\section{Introduction}

A recent study \cite{SSDIF} characterized the intersection of two principal
classes of one-dimensional random motions: the class of selfsimilar
processes, and the class of Ito diffusion processes. The former class
manifests random motions whose trajectories are, statistically, fractal
objects \cite{EM}. The latter class manifests random motions whose dynamics
are governed by Ito's stochastic differential equation \cite{Fri}-\cite{Arn}%
. Both these classes are of major importance, theoretical and practical
alike. In this paper, progressing from the theoretical study \cite{SSDIF} to statistical-physics applications, we explore a practical model for \emph{anomalous
diffusion}. 

Regular and anomalous diffusion assume central roles in science and
engineering \cite{PSW}-\cite{Shl}. The paradigmatic model for regular
diffusion is Brownian motion \cite{Gar}-\cite{Kam}, which is a mathematical
object of exquisite beauty and riches \cite{KaS}. Maintaining key properties
of Brownian motion -- namely, being a selfsimilar process with finite variance and
with a continuous trajectory -- the paradigmatic model for anomalous
diffusion is \emph{fractional Brownian motion} \cite{CheG}-\cite{CWMS}. On
the flip side, fractional Brownian motion fails to maintain other key
properties of Brownian motion; specifically, fractional Brownian motion is
neither a Markov process, nor a martingale \cite{EM}. Consequently, both the
numerical simulation and the analytic tracking of fractional Brownian motion
are challenging tasks.

In analogy with fractional Brownian motion, we term this paper's anomalous-diffusion model \emph{fractional Ito motion}. As fractional
Brownian motion, fractional Ito motion also maintains the aforementioned key
properties of Brownian motion: it is a selfsimilar process with finite
variance and with a continuous trajectory. In sharp contrast to fractional
Brownian motion, fractional Ito motion is a Markov process, as well as a
martingale. Consequently, fractional Ito motion is a highly useful and a
highly workable model for anomalous diffusion.

A Stratonovich counterpart of fractional Ito motion -- in which the underpinning Ito integration is replaced by Stratonovich integration -- was explored in \cite{CCM}. A special case of fractional Ito motion was investigated in \cite{LB}, and, most recently, fractional Ito motion was investigated in the context of stochastic resetting \cite{LLG}. In this paper, taking on a selfsimilarity approach, we conduct a detailed
comparison between two anomalous diffusion models: \emph{fractional Brownian motion} vs. \emph{fractional Ito motion}. 

The paper is organized as follows. We begin with an overview of regular and
anomalous diffusion, and with a concise description of the three aforementioned
models (section \ref{1}): Brownian motion, fractional Brownian motion, and fractional Ito
motion. We then present the detailed comparison between fractional Brownian motion and fractional Ito motion (section \ref{2}). Thereafter, we present a ``diffusion-in-a-logarithmic-potential'' representation of fractional Ito motion (section \ref{3}), and conclude with a summary (section \ref{4}).

\section{\label{1}Regular and anomalous diffusion}

Diffusion is a generic name for random motions that diffuse with time.
Arguably, the most common method to measure diffusivity is mean square
displacement (MSD). To describe this method, consider a one-dimensional
random motion whose trajectory is $X\left( t\right) $ ($t\geq 0$). Namely,
the motion initiates at time $0$, and its position at time $t$ is $X\left(
t\right) $ (a point on the real line). Hence, at time $t$, the motion's
displacement -- relative to its initial position -- is $\left\vert X\left(
t\right) -X\left( 0\right) \right\vert $. In turn, the motion's MSD at time $%
t$ is $\phi \left( t\right) =\mathbf{E}[\left\vert X\left( t\right) -X\left(
0\right) \right\vert ^{2}]$, where here and hereinafter $\mathbf{E}\left[
\cdot \right] $ denotes the operation of statistical expectation. Rather
generally, the random motion under consideration can be said to be a \emph{%
diffusion} if: its MSD function $\phi \left( t\right) $ increases with time,
and diverges as time grows infinitely large.

Often, diffusions take place at microscopic spatio-temporal scales, but yet
they are observed at macroscopic spatio-temporal scales. As explained in the
Methods, on the macroscopic level the only admissible form of the MSD
function is a power-law: 
\begin{equation}
\phi \left( t\right) =c\cdot t^{\epsilon },  \label{10}
\end{equation}%
where $c$ is a positive coefficient, and where $\epsilon $ is a positive
exponent. A broad and rich framework for random motions that yield such a
power-law MSD form are \emph{selfsimilar processes} with finite variance 
\cite{EM}. To describe the notion of selfsimilarity consider, as above, a
one-dimensional random motion whose trajectory is $X\left( t\right) $ ($%
t\geq 0$).

The random motion under consideration is selfsimilar if -- for any given
positive scale $s$ -- its two following versions are statistically equal 
\cite{EM}: the time-scaled version $X\left( st\right) $ ($t\geq 0$); and the
space-scaled version $s^{H}X\left( t\right) $ ($t\geq 0$), where $H$ is a
positive \emph{Hurst exponent}. Namely, with regard to the random motion
under consideration, selfsimilarity means that: changing the underlying
temporal scale from $1$ to $s$ is statistically equivalent to changing the
underlying spatial scale from $1$ to $s^{H}$. Selfsimilarity implies that
the motion initiates at the spatial origin $X\left( 0\right) =0$, and that
the random variable $X\left( t\right) $ is equal in law to the random
variable $t^{H}\cdot X\left( 1\right) $. Consequently -- in the
finite-variance case -- selfsimilarity yields the power-law MSD form of Eq. (%
\ref{10}) with: coefficient $c=\mathbf{E}[\left\vert X\left( 1\right)
\right\vert ^{2}]$, and exponent $\epsilon =2H$.

The above discussion is mathematically neat, yet it is merely theoretical.
Does the real world corroborate the theory? The answer is a resounding yes:
real-world observations happen to be in perfect accord with the theory.
Indeed, pioneered by the experimental observations of Jan Ingen-Housz \cite%
{Vpas}, of Sir Robert Brown \cite{Bro}, and of Jean Perrin \cite{Per}, a
vast body of empirical evidence established the ubiquity of diffusions with
linear MSD functions, $\phi \left( t\right) =c\cdot t$. Such diffusions are
so prevalent that they were assumed to manifest diffusive motions at large.
In turn, the convention was to quantify a given diffusive motion by its
\textquotedblleft diffusion coefficient\textquotedblright\ -- the slope of
its linear MSD function.

In the nineteen-seventies of the twentieth century new experimental
observations began to challenge the then commonplace assumption that
diffusive motions have linear MSD functions \cite{SL1}-\cite{SMon}. As the
new body of empirical evidence grew larger and larger, it became clear that
in addition to the well-known diffusive motions with linear MSD functions,
there are also diffusive motions with power-law MSD functions. Consequently,
diffusions were distinguished by three different categories: \emph{regular
diffusion}, characterized by linear MSD functions; \emph{sub-diffusion},
characterized by power-law MSD functions with sub-linear exponents $\epsilon
<1$; and \emph{super-diffusion}, characterized by power-law MSD functions
with super-linear exponents $\epsilon >1$.

The experimental discovery of sub-diffusive and super-diffusive random
motions ushered in the multidisciplinary scientific field of \emph{anomalous
diffusion} \cite{PSW}-\cite{Shl}. This field attracted substantial
scientific interest \cite{GAA}-\cite{Metz}, and the scientific exploration
of the field led to the conclusions that \textquotedblleft anomalous is
normal\textquotedblright\ \cite{SLLS} and that \textquotedblleft anomalous
is ubiquitous\textquotedblright\ \cite{AnoUbi}.

Over forty years after its inception, anomalous diffusion continues to be a
vibrant multidisciplinary scientific field. Examples of recent
anomalous-diffusion research include: anomalous diffusion in complex media 
\cite{VSSP}; anomalous diffusion in the evolution of inhomogeneous systems 
\cite{OFL}; anomalous diffusion in random dynamical systems \cite{SKla};
anomalous diffusion and recurrent neural networks \cite{BSFV}; heterogeneous
diffusion processes \cite{WCLM}; anomalous diffusion in heterogeneous binary
media \cite{FTPB}; anomalous diffusion and time-dependent diffusivity \cite%
{CSRM}; anomalous diffusion in the noisy voter model \cite{Kon}-\cite{KazK};
anomalous diffusion under stochastic resetting \cite{MasM}-\cite{FauS};
anomalous diffusion in comb structures \cite{DzS}-\cite{Iom}; and random
diffusivity \cite{WSSC}-\cite{SanM}.

We now turn to describe the three random-motion models that were noted in
the introduction. We begin with the paradigmatic models of Brownian motion
and fractional Brownian motion. Thereafter, we present the model of
fractional Ito motion.

\subsection{\label{1A}Brownian motion and fractional Brownian motion}

Following the trailblazing theoretical works of Louis Bachelier in finance 
\cite{Bac1}-\cite{Bac2}, of Albert Einstein and Marian Smoluchowski in
physics \cite{Ein}-\cite{Smo}, and of Norbert Wiener in mathematics \cite%
{Wie}, a particular random motion emerged as the paradigmatic model for
regular diffusion in science and engineering \cite{Gar}-\cite{Kam}: \emph{%
Brownian motion} (BM), $B\left( t\right) $ ($t\geq 0$). Named in honor of
Sir Robert Brown, and also termed \emph{Wiener process} in honor of Norbert
Wiener, BM is a profound object that exhibits a host of mathematical,
statistical, and geometric properties \cite{KaS}.

In this paper we shall address the following key BM\ properties. I) BM has
finite variance, i.e. its positions have finite first-order and second-order
moments; hence BM has a well-defined MSD function. II) BM is a symmetric
process, i.e. its trajectory is statistically identical to its mirror
trajectory, $-B\left( t\right) $ ($t\geq 0$). III) The trajectory of BM\ is
continuous. IV) BM is a selfsimilar process with Hurst exponent $H=\frac{1}{2%
}$; hence the MSD function of BM is linear, $\epsilon=1$, and hence it is a regular
diffusion indeed. V) BM is a Gaussian process, i.e. its finite-dimensional
distributions are multivariate Normal. VI) BM is a Levy process, i.e. its
increments are stationary and its non-overlapping increments are
independent. VII) BM is a Markov process, i.e. -- at any given time point $t$
-- its future trajectory depends only at its present position, $B\left(
t\right) $, and it does not depend on its past trajectory. VIII) BM is a
martingale.

The properties of BM\ imply that its positions, as well as its increments,
are Normal random variables with zero means. The inherent scale of BM is set
so that the random variable $B\left( 1\right) $ -- the position of BM at the
time point $1$ -- is \textquotedblleft standard Normal\textquotedblright ,
i.e.: $B\left( 1\right) $ is a Normal random variable with zero mean and
with unit variance.

In order to model anomalous diffusion one needs to go beyond BM. To do so --
within the realm of symmetric and selfsimilar processes with finite variance
and continuous trajectories -- leads to the following model: \emph{%
fractional Brownian motion} (FBM), $B_{H}\left( t\right) $ ($t\geq 0$),
where the subscript $H$ manifests the underlying Hurst exponent. Pioneered
by Andrey Kolmogorov \cite{Kol}, by Akiva Yaglom \cite{Yag}, and by Benoit
Mandelbrot and John an Ness \cite{MV}, FBM is a well established
generalization of BM \cite{Mish}-\cite{BMRS}. The Hurst exponent of FBM\
takes values in the range $0<H<1$, and hence the diffusivity of FBM\ is as
follows: sub-diffusion in the exponent range $0<H<\frac{1}{2}$;
super-diffusion in the exponent range $\frac{1}{2}<H<1$; and regular
diffusion at the exponent value $H=\frac{1}{2}$ -- in which case FBM is BM.

One the one hand -- as BM -- FBM is a Gaussian process, its trajectory is
continuous, and its increments are stationary \cite{EM}. On the other hand
-- in sharp contrast to BM -- the non-overlapping increments of FBM\ are
dependent, and FBM\ is neither a Markov process nor is it a martingale \cite%
{EM}. Moreover, FBM is not even a semi-martingale, and this fact has major
implications in the context of stochastic integration \cite{EM}.

To gain insight regarding the difference between BM\ and FBM, we turn to
their velocities. The velocity of BM, $\dot{B}\left( t\right) $ ($t\geq 0$),
is commonly known as \emph{white noise}. The velocity of FBM, $\dot{B}%
_{H}\left( t\right) $ ($t\geq 0$), is a moving average -- with a power-law
kernel -- of white noise. Specifically, the velocity of FBM with Hurst
exponent $H\neq \frac{1}{2}$ admits the following moving-average
representation \cite{FracMo}:%
\begin{equation}
\dot{B}_{H}\left( t\right) =\left( H-\frac{1}{2}\right) \int_{-\infty
}^{t}\left( t-u\right) ^{H-\frac{3}{2}}\dot{B}\left( u\right) du.
\label{11}
\end{equation}%
The moving-average representation of Eq. (\ref{11}) actually uses two
independent white noises: one defined over the negative time axis, $\dot{B}%
\left( u\right) $ ($u<0$); and one defined over the non-negative time axis, $%
\dot{B}\left( u\right) $ ($u\geq 0$).

As FBM is a selfsimilar process, it initiates at the spatial origin $%
B_{H}\left( 0\right) =0$. Hence, integrating Eq. (\ref{11}) yields%
\begin{equation}
B_{H}\left( t\right) =\int_{-\infty }^{0}\left[ \left( t-u\right) ^{H-\frac{1%
}{2}}-\left( 0-u\right) ^{H-\frac{1}{2}}\right] \dot{B}\left( u\right)
du+\int_{0}^{t}\left( t-u\right) ^{H-\frac{1}{2}}\dot{B}\left( u\right) du%
.  \label{12}
\end{equation}%
The integral representation of Eq. (\ref{12}) comprises two parts: an
integral of the white noise $\dot{B}\left( u\right) $ over the negative
temporal axis, $-\infty <u<0$; and an integral of the white noise $\dot{B}%
\left( u\right) $ over the temporal interval $0\leq u\leq t$. The integral
representation of Eq. (\ref{12}) covers also the Hurst exponent $H=\frac{1}{2%
}$, and setting this Hurst exponent in Eq. (\ref{12}) yields BM: $%
B_{1/2}\left( t\right) =B\left( t\right) $. Hence, FBM is indeed a
generalization of BM. Schematic illustrations of the integration kernel of
Eq. (\ref{12}) -- whose shape is determined by the value of the Hurst
exponent $H$ -- are depicted in Figure 1.

As noted above, FBM is not a Markov process. There are one-dimensional
random motions that are non-Markov, but yet they can be transformed to a
Markov-process representation by elevating from a single dimension to
several dimensions. FBM\ is a highly non-Markov process in the following
sense: in order to transform FBM\ to a Markov-process representation one has
to elevate from a single dimension to infinitely many dimensions \cite%
{FracMo}. Consequently, the numerical simulation of FBM\ is a challenging
computational task \cite{Yin}-\cite{SMB}.

So, on the one hand, FBM is a random-motion model -- which is selfsimilar
and Gaussian, whose positions have a finite variance, and whose trajectory
is continuous -- that produces both sub-diffusion and super-diffusion. On
the other hand, working with FBM is quite challenging from the perspectives
of stochastic integration and numerical simulation, as well as from the
perspective of analytic tractability. The intricacy of FBM\ gives rise to
the following question: is there an alternative anomalous-diffusion model
that shares the `upside' features of FBM, but not the `downside' features of
FBM? The answer, as we shall now argue, is affirmative indeed.

\subsection{\label{1B}Fractional Ito motion}

Consider a one-dimensional deterministic motion $X\left( t\right) $ ($t\geq
0 $), whose dynamics are governed by the ordinary differential equation
(ODE) $\dot{X}\left( t\right) =\mu \left[ X\left( t\right) \right] $.
Namely, when the deterministic motion is at the position $x$ (a point on the
real line) then its velocity is $\mu \left( x\right) $ (a real number).
Adding white noise to the dynamics changes the motion from deterministic to
random, and changes the ODE to a Langevin stochastic differential equation
(SDE) \cite{Lan}-\cite{Pav}: $\dot{X}\left( t\right) =\mu \left[ X\left(
t\right) \right] +\nu \dot{B}\left( t\right) $, where $\nu $ is a positive
parameter that manifest the white-noise magnitude. In turn, replacing the
constant white-noise magnitude $\nu $ by a position-dependent magnitude
changes the Langevin SDE to an Ito SDE \cite{Fri}-\cite{Arn},\cite{Ito}-\cite%
{IK}: $\dot{X}\left( t\right) =\mu \left[ X\left( t\right) \right] +\sigma
\lbrack X\left( t\right) ]\dot{B}\left( t\right) $. In the jargon of
mathematical finance \cite{Hul}: when the random motion is at the position $%
x $ (a point on the real line) then its \textquotedblleft
drift\textquotedblright\ is $\mu \left( x\right) $ (a real number), and its
\textquotedblleft volatility\textquotedblright\ is $\sigma \left( x\right) $
(a positive number).

With the notion of Ito SDEs recalled, we introduce the following model: 
\emph{fractional Ito motion} (FIM), $I_{H}\left( t\right) $ ($t\geq 0$),
where the subscript $H$ manifests the underlying Hurst exponent. In analogy
with FBM, we name FIM in honor of Kioshi Ito -- the mathematician that
invented the white-noise stochastic calculus. The dynamics of FIM are
governed by the Ito SDE%
\begin{equation}
\dot{I}_{H}\left( t\right) =\left\vert I_{H}\left( t\right) \right\vert ^{1-
\frac{1}{2H}}\dot{B}\left( t\right)   \label{21}
\end{equation}%
Namely, FIM has a zero drift, $\mu \left( x\right) =0$, and a power-law
volatility, $\sigma \left( x\right) =\left\vert x\right\vert ^{1-\frac{1}{2H}%
}$. In general, a random motion whose dynamics are governed by an Ito SDE is
a Markov process, and its trajectory is continuous \cite{KaS}; hence, in
particular, FIM exhibits these properties.

We initiate FIM from the spatial origin $I_{H}\left( 0\right) =0$. Hence,
integrating Eq. (\ref{21}) yields%
\begin{equation}
I_{H}\left( t\right) =\int_{0}^{t}\left\vert I_{H}\left( u\right)
\right\vert ^{1-\frac{1}{2H}}\dot{B}\left( u\right) du.  \label{22}
\end{equation}%
The right-hand side of Eq. (\ref{22}) is a running Ito integral. A general
running Ito integral -- with an integrand that does not `look into the
future' -- is a symmetric process and a martingale \cite{KaS}; hence, in
particular, FIM exhibits these properties.

According to general results regarding selfsimilar diffusions \cite{SSDIF},
FIM is a selfsimilar process, and the Hurst exponent of FIM takes values in
the range $0<H<1$. Hence, the diffusivity of FIM\ is identical to the
aforementioned diffusivity of FBM: sub-diffusion in the exponent range $0<H<%
\frac{1}{2}$; super-diffusion in the exponent range $\frac{1}{2}<H<1$; and
regular diffusion at the exponent value $H=\frac{1}{2}$ -- in which case FIM
is BM. Indeed, setting $H=\frac{1}{2}$ in Eq. (\ref{22}) yields BM: $%
I_{1/2}\left( t\right) =B\left( t\right) $. Schematic illustrations of the
FIM `volatility landscape'\ -- whose shape is determined by the value of the
Hurst exponent $H$ -- are depicted in Figure 1.

So, on the one hand, FIM shares the `upside' features of FBM: it is a
random-motion model -- which is symmetric and selfsimilar, and whose
trajectory is continuous -- that generalizes BM, and that produces both
sub-diffusion and super-diffusion. On the other hand, as FIM is a Markov
process and a martingale, it circumvents the `downside' features of FBM: it
well applies in the context of stochastic integration, its numerical
simulation is easy and straightforward, and it is analytically tractable.

Replacing the Ito SDE (\ref{21}) by an identical Stratonovich SDE results in a ``fractional Stratonovich motion'' counterpart of FIM; this counterpart of the FIM model was explored in \cite{CCM}. A special case of the FIM model was investigated in \cite{LB}; in this special case the motion runs over the positive half-line ($0<x<\infty $), and the Hurst exponent is in the range $\frac{1}{2} \leq H <1$. Most recently, FIM was also investigated in the context of stochastic resetting \cite{LLG}. To the best of our knowledge, the detailed FBM vs. FIM comparison -- which is conducted and presented in this paper -- is entirely new.  

\begin{figure}[H]
\begin{centering}
\includegraphics[width=0.5\paperwidth,height=0.5\paperwidth,keepaspectratio]{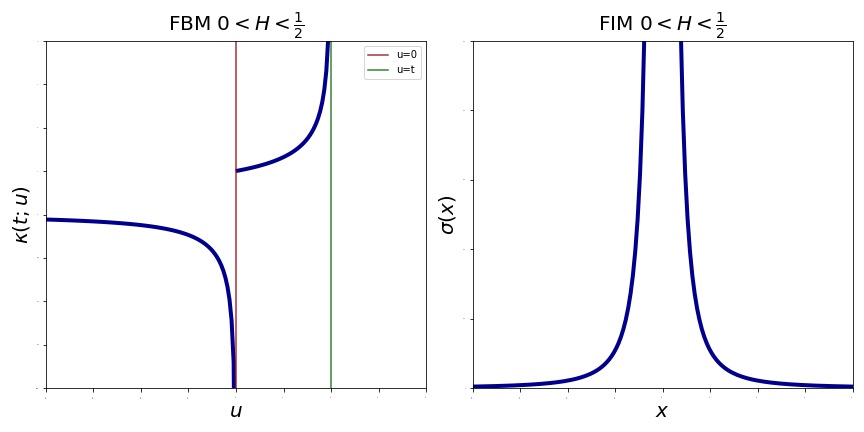}
\includegraphics[width=0.5\paperwidth,height=0.5\paperwidth,keepaspectratio]{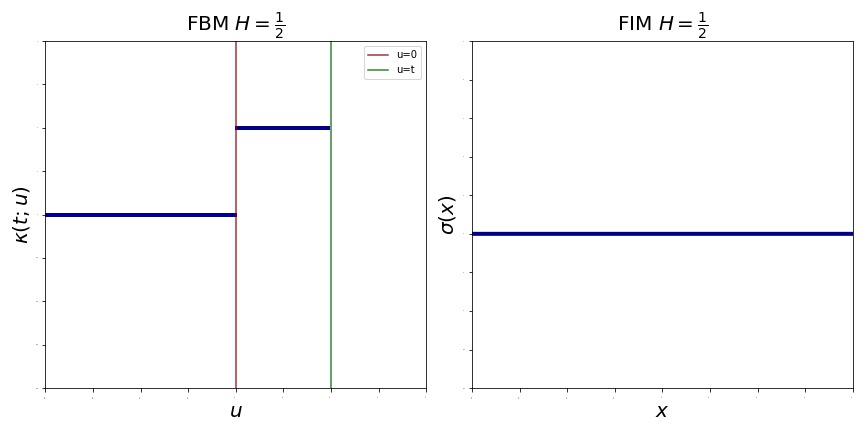}
\includegraphics[width=0.5\paperwidth,height=0.5\paperwidth,keepaspectratio]{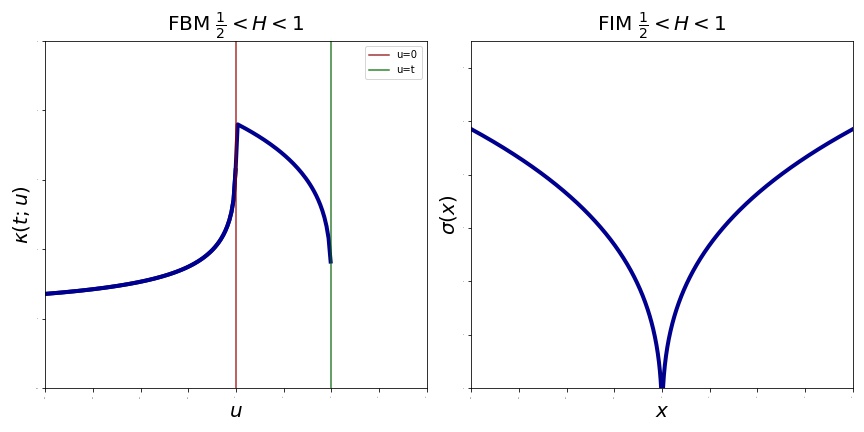}
\par\end{centering}
\caption{schematic illustrations of the shapes of the FBM integration kernel, and of the FIM
volatility landscape. Depicted in the left panels, the FBM integration kernel (for a fixed positive time point $t$) is: $\kappa \left( t;u\right) =\left( t-u\right) ^{H-\frac{1}{2}}-\left( 0-u\right) ^{H-\frac{1}{2}}$ over the negative temporal axis, $%
-\infty <u<0$; and $\kappa \left( t;u\right) =\left( t-u\right) ^{H-\frac{1}{%
2}}$ over the temporal interval $0\leq u\leq t$.  Depicted in the right panels, the FIM volatility
landscape is: $\sigma \left( x\right) =\left\vert x\right\vert ^{1-\frac{1}{2H%
}}$ over the spatial axis $-\infty <x<\infty $. The temporal
kernel function $\kappa \left( t;u\right) $ and the spatial volatility
function $\sigma \left( x\right) $ display markedly different shapes in the
three different diffusivity categories: sub-diffusion $0<H<\frac{1}{2}$ (top panels); super-diffusion $\frac{1}{2}<H<1$ (bottom panels); and regular diffusion $H=\frac{1}{2}$ (middle panels).}
\end{figure}

\section{\label{2}Fractional Brownian motion vs. fractional Ito motion}

Describing the FBM\ model and FIM model, the previous section revealed some
of the marked differences between these two anomalous-diffusion models. In
this section we carry on with a series of analytic comparisons that examine
and pinpoint additional profound differences between FBM and FIM. Visual comparisons between simulated trajectories of FBM and FIM are offered by Figure 2 (for sub-diffusion) and by Figure 3 (for super-diffusion).

\ \ 

\begin{figure}[H]
\begin{centering}

\includegraphics[width=0.7\paperwidth,height=0.7\paperwidth,keepaspectratio]{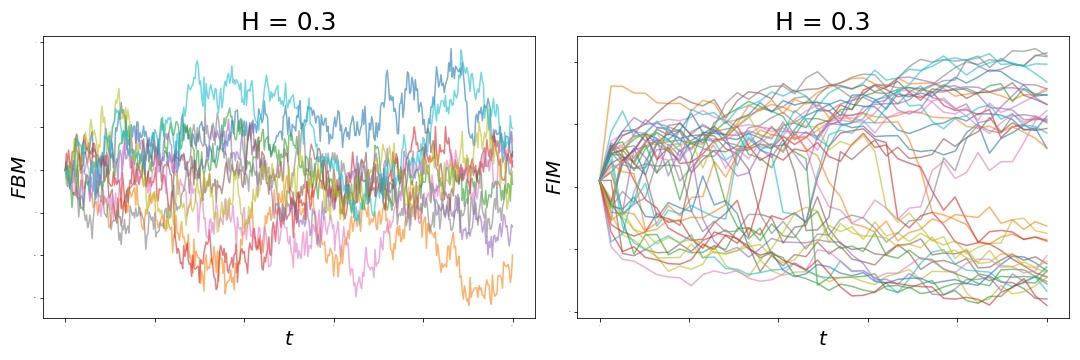}
\includegraphics[width=0.7\paperwidth,height=0.7\paperwidth,keepaspectratio]{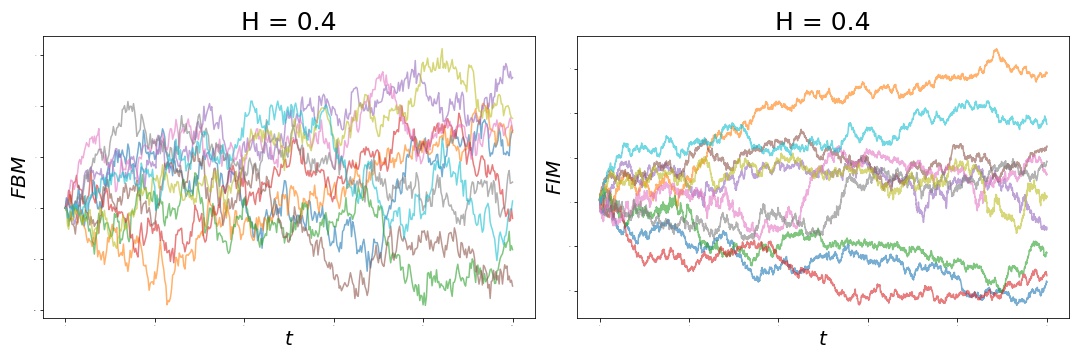}
\par\end{centering}
\caption{simulated trajectories of sub-diffusive FBM (left panels) and of sub-diffusive FIM (right panels).}

\end{figure}

\begin{figure}[H]
\begin{centering}
\includegraphics[width=0.7\paperwidth,height=0.7\paperwidth,keepaspectratio]{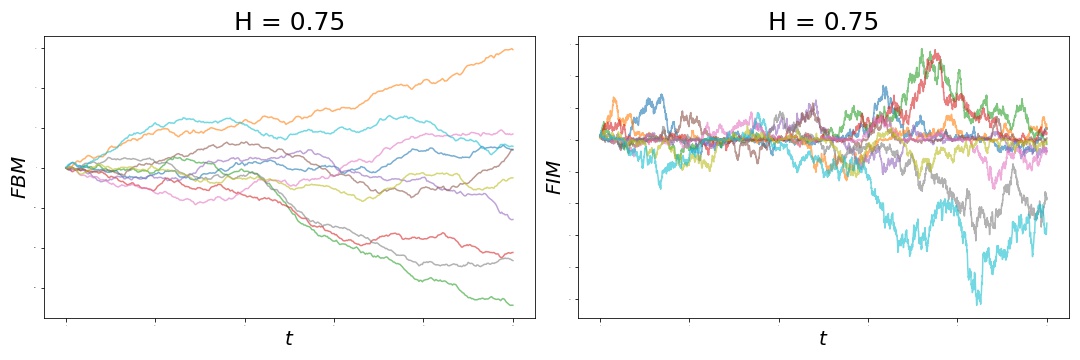}
\includegraphics[width=0.7\paperwidth,height=0.7\paperwidth,keepaspectratio]{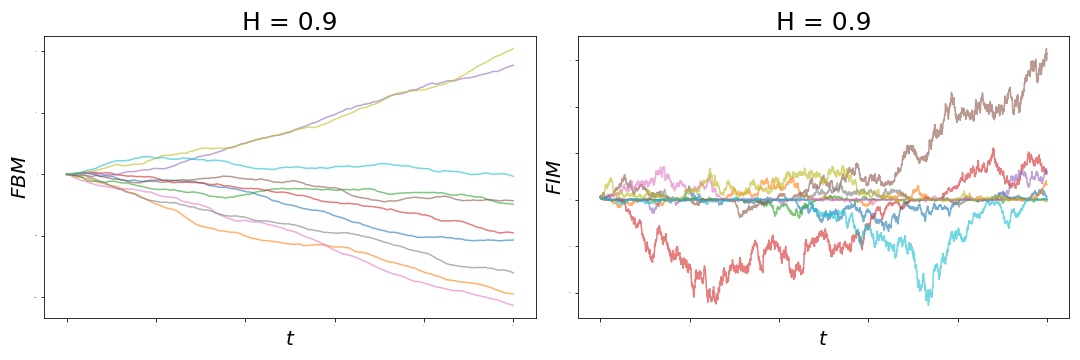}
\par\end{centering}
\caption{simulated trajectories of super-diffusive FBM (left panels) and of super-diffusive FIM (right panels).}

\end{figure}

\subsection{\label{A}Dissipation comparison}

As FBM\ and FIM are selfsimilar processes, they both initiate from the
spatial origin. Hence, from a probabilistic perspective, at time $0$ both
these random motions manifest a unit mass that is placed at the origin. In
turn, at a positive time $t$, this unit mass dissipates, and the `shape of
the dissipation' is quantified by a probability density function: the
density of the random variable $B_{H}\left( t\right) $, in the case of FBM;
and the density of the random variable $I_{H}\left( t\right) $, in the case
of FIM.

The properties of FBM imply that the random variable $B_{H}\left( t\right) $
is Normal with mean zero with variance $\mathbf{Var}\left[ B_{H}\left(
t\right) \right] =\mathbf{Var}\left[ B_{H}\left( 1\right) \right] \cdot
t^{2H}$. Hence, setting $b=\mathbf{Var}\left[ B_{H}\left( 1\right) \right] $%
, the density of the random variable $B_{H}\left( t\right) $ is%
\begin{equation}
\frac{1}{\sqrt{2\pi b}}\cdot \frac{1}{t^{H}}\exp \left( -\frac{x^{2}}{%
2bt^{2H}}\right)  \label{31}
\end{equation}%
($-\infty <x<\infty $). This density is a symmetric `bell curve' (see Figure
4): it vanishes at $x\rightarrow \pm \infty $, and it has a unimodal shape
the peaks at the spatial origin. We emphasize that the shape of this
density is the same for all the values of the Hurst exponent $H$.

It follows from a general result established in \cite{SSDIF} that the
density of the random variable $I_{H}\left( t\right) $ is%
\begin{equation}
\frac{1}{2H\Gamma \left( 1-H\right) }\cdot \left( \frac{2H^{2}}{t}\right)
^{1-H}\exp \left( -\frac{2H^{2}}{t}\left\vert x\right\vert ^{\frac{1}{H}%
}\right) \left\vert x\right\vert ^{\frac{1}{H}-2}  \label{32}
\end{equation}%
($-\infty <x<\infty $). This density is symmetric, and it vanishes at $%
x\rightarrow \pm \infty $. The shape of this density is determined by the
value of the Hurst exponent $H$, as follows (see Figure 4).

\begin{enumerate}
\item[$\bullet $] In the sub-diffusion range, $0<H<\frac{1}{2}$, the density
has a bimodal shape: it vanishes at the spatial origin, and it peaks at the
spatial points $\pm \left( \frac{1-2H}{2H^{2}}\right) t^{H}$.

\item[$\bullet $] At the regular-diffusion value, $H=\frac{1}{2}$, the
density is a `bell curve': it has a unimodal shape the peaks at the spatial
origin.

\item[$\bullet $] In the super-diffusion range, $\frac{1}{2}<H<1$, the
density has a unimodal shape the explodes at the spatial origin.
\end{enumerate}

Evidently, the differences between the shape of the Gaussian FBM density of
Eq. (\ref{31}) and the shape of the non-Gaussian (for $H\neq \frac{1}{2}$)
FIM density of Eq. (\ref{32}) are dramatic. On the one hand, changing the
Hurst exponent $H$ in the FBM model has no qualitative effect on the shape
of the dissipation pattern. On the other hand, changing the Hurst exponent $%
H $ in the FIM model has a profound qualitative effect on the shape of the
dissipation pattern. 

The dissipation pattern of ``fractional Stratonovich motion'' -- the Stratonovich counterpart of FIM -- displays a behavior which is similar to that of the FIM density of Eq. (\ref{32}) \cite{CCM}. Non-Gaussian diffusions (i.e. diffusions with non-Gaussian dissipation patterns), such as FIM and its Stratonovich counterpart, attracted substantial scientific interest recently \cite{CSM}-\cite{Met2}.

\begin{figure}[H]
\begin{centering}
\includegraphics[width=0.5\paperwidth,height=0.5\paperwidth,keepaspectratio]{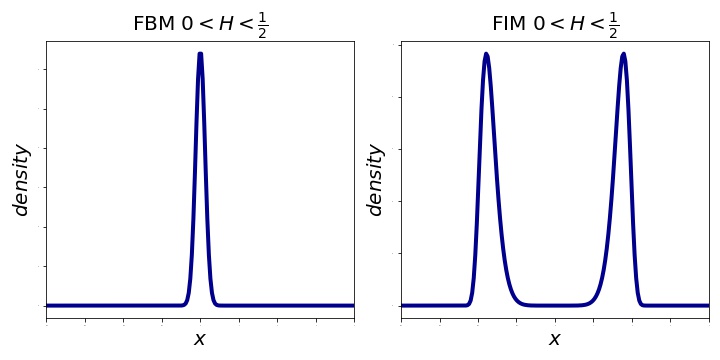}
\includegraphics[width=0.5\paperwidth,height=0.5\paperwidth,keepaspectratio]{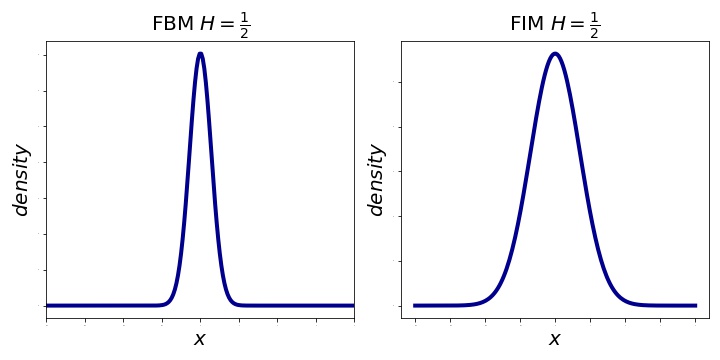}
\includegraphics[width=0.5\paperwidth,height=0.5\paperwidth,keepaspectratio]{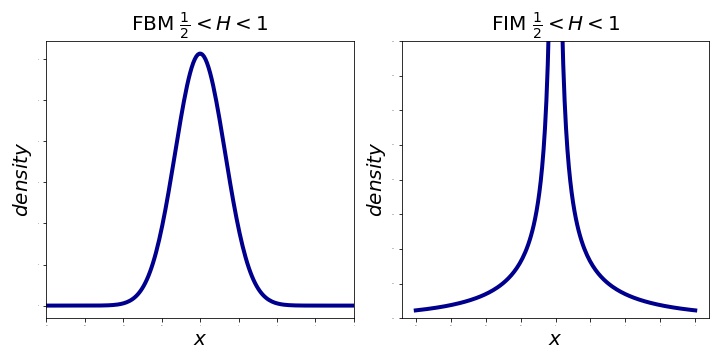}
\par\end{centering}
\caption{schematic illustrations of the shapes of the FBM and FIM dissipation patterns (for a fixed positive time point $t$). Depicted in the left panels, the FBM dissipation pattern is the probability density function of the random variable $B_{H}\left( t\right)$ (Eq. (\ref{31})). Depicted in the right panels, the FIM dissipation pattern is the probability density function of the random variable $I_{H}\left( t\right)$ (Eq. (\ref{32})). While the FBM dissipation pattern displays the same unimodal `bell-curve' shape, the FIM dissipation pattern displays markedly different shapes in the three different diffusivity categories: bimodal for sub-diffusion $0<H<\frac{1}{2}$ (top panels); unimodal and explosive for super-diffusion $\frac{1}{2}<H<1$ (bottom panels); and unimodal `bell-curve' for regular diffusion $H=\frac{1}{2}$ (middle panels).}

\end{figure}

\begin{figure}[H]
\begin{centering}
\includegraphics[width=0.5\paperwidth,height=0.5\paperwidth,keepaspectratio]{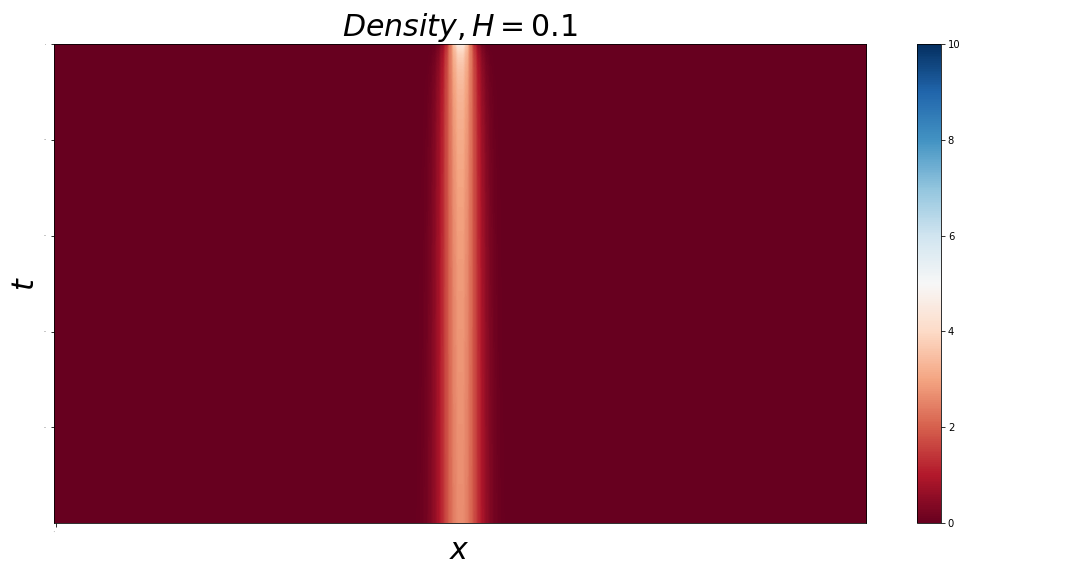}
\includegraphics[width=0.5\paperwidth,height=0.5\paperwidth,keepaspectratio]{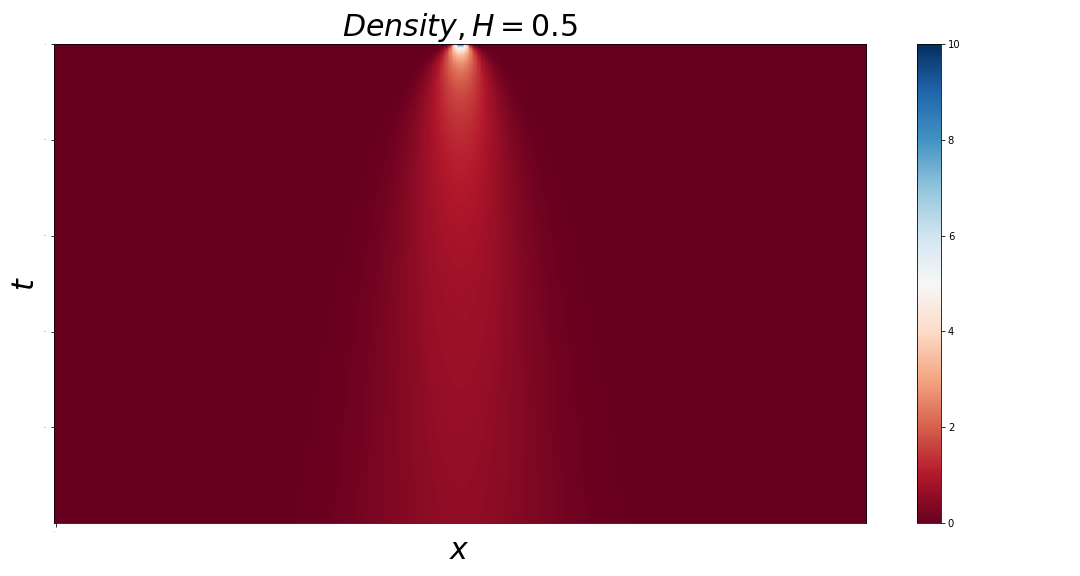}
\includegraphics[width=0.5\paperwidth,height=0.5\paperwidth,keepaspectratio]{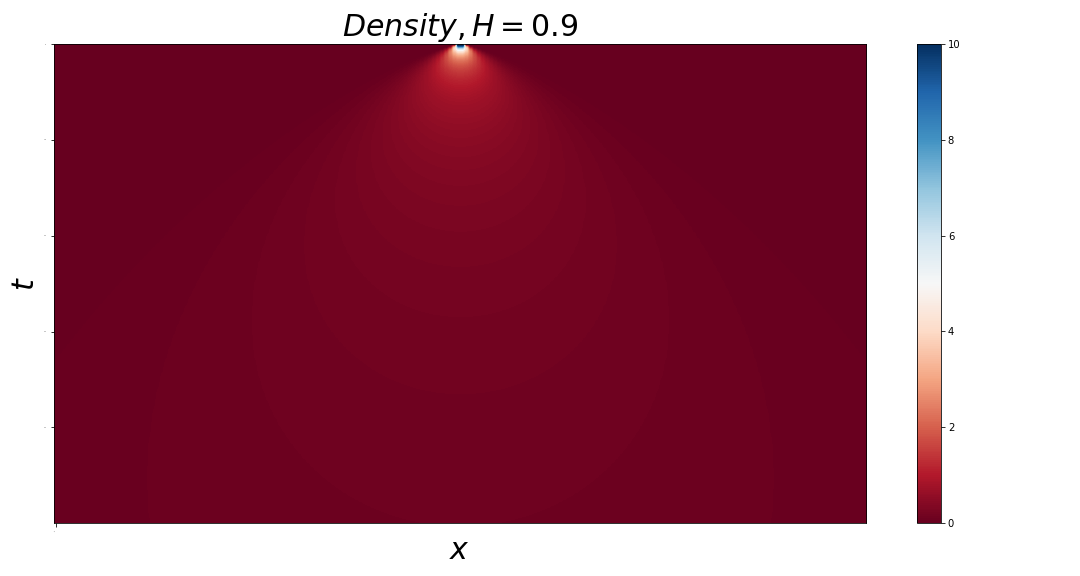}
\par\end{centering}
\caption{`heat maps' for three examples of the FBM dissipation pattern -- the probability density function of the random variable $B_{H}\left( t\right)$ (Eq. (\ref{31})). The examples correspond to the three different diffusivity categories of FBM: sub-diffusion (top panel); regular diffusion (middle panel); super diffusion (bottom panel).}

\end{figure}

\begin{figure}[H]
\begin{centering}
\includegraphics[width=0.5\paperwidth,height=0.5\paperwidth,keepaspectratio]{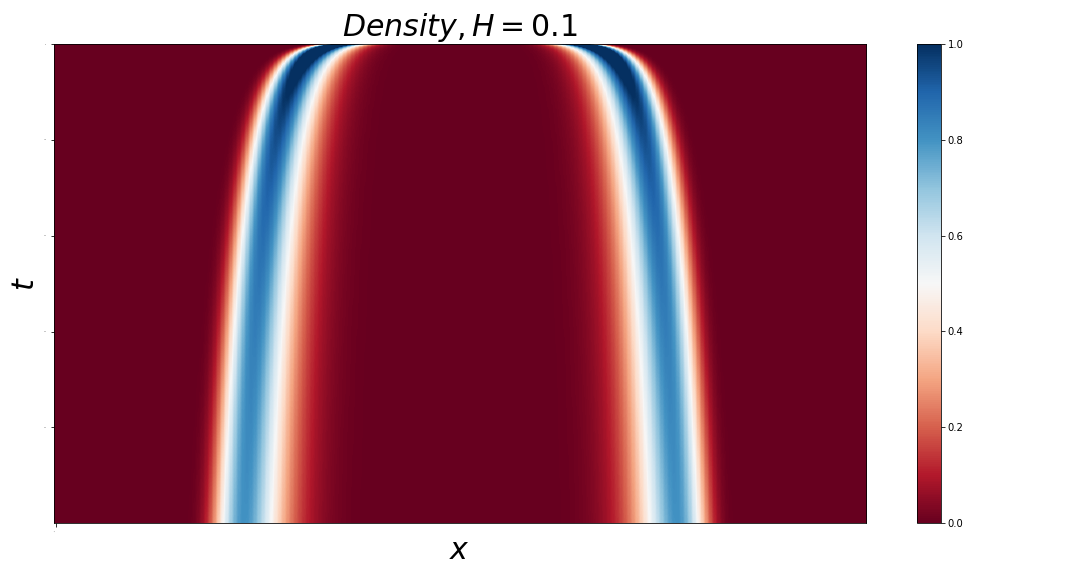}
\includegraphics[width=0.5\paperwidth,height=0.5\paperwidth,keepaspectratio]{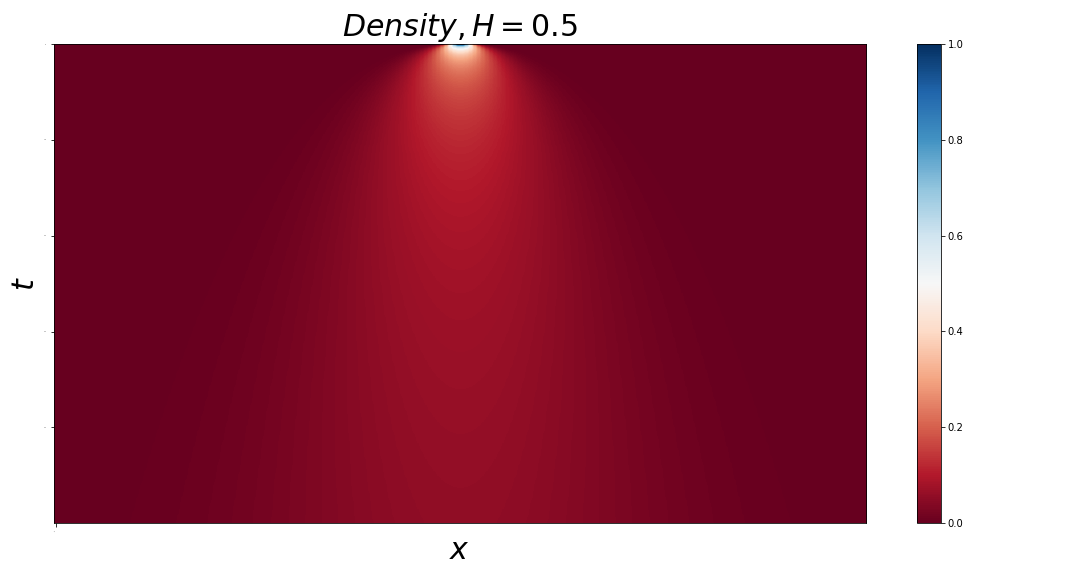}
\includegraphics[width=0.5\paperwidth,height=0.5\paperwidth,keepaspectratio]{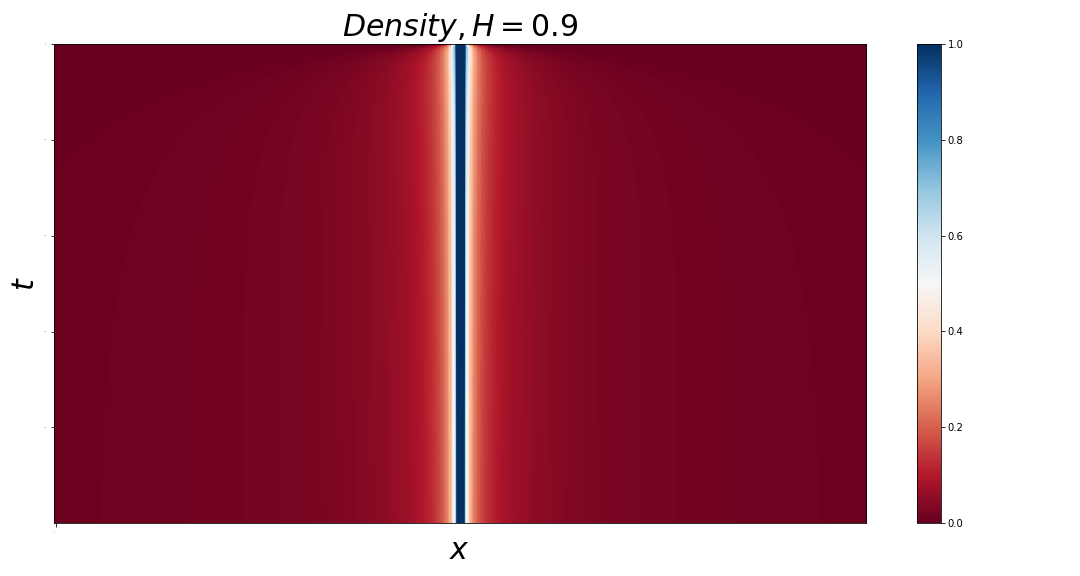}
\par\end{centering}
\caption{`heat maps' for three examples of the FIM dissipation pattern -- the probability density function of the random variable $I_{H}\left( t\right)$ (Eq. (\ref{32})). The examples correspond to the three different diffusivity categories of FIM: sub-diffusion (top panel); regular diffusion (middle panel); super diffusion (bottom panel).}

\end{figure}

\subsection{\label{B}Tail-behavior comparison}

Consider a symmetric one-dimensional random motion whose trajectory is $%
X\left( t\right) $ ($t\geq 0$). The tail behavior of the random motion under
consideration is the asymptotic behavior of the probability $\Pr \left[
\left\vert X\left( t\right) \right\vert >l\right] $ in the limit $%
l\rightarrow \infty $. Namely, the tail behavior quantifies the asymptotic
likelihood of the following rare event: the motion's displacement, at time $%
t $, relative to its initial position, is greater than the level $l>>1$.

Using Eqs. (\ref{31})-(\ref{32}) and L'Hospital's rule, we conduct three
tail-behavior comparisons: FBM vs. BM; FIM vs. BM; and FIM vs. FBM. These
tail-behavior comparisons are presented in Table 1, and they yield the
following `asymptotic orderings'. In the case of sub-diffusion ($0<H<\frac{1%
}{2}$): the tail behavior of FIM is infinitely `lighter' than the tail
behavior of FBM -- which, in turn, is infinitely `lighter' than the tail
behavior of BM. In the case of super-diffusion ($\frac{1}{2}<H<1$): the tail
behavior of FIM is infinitely `heavier' than the tail behavior of FBM --
which, in turn, is infinitely `heavier' than the tail behavior of BM.

\ \ \ \ 

\begin{center}
{\large Table 1}
\end{center}

\ \ \ \ \ 

\begin{tabular}{|c|c|c|c|}
\hline
\textbf{Comparison} & \textbf{Limit} & $%
\begin{array}{c}
\text{\textbf{Sub diffusion}} \\ 
\text{(}0<H<\frac{1}{2}\text{)}%
\end{array}%
$ & $%
\begin{array}{c}
\text{\textbf{Super diffusion}} \\ 
\text{(}\frac{1}{2}<H<1\text{)}%
\end{array}%
$ \\ \hline
1) FBM vs. BM & $%
\begin{array}{c}
\text{ } \\ 
\lim_{l\rightarrow \infty }\frac{\Pr \left[ \left\vert B_{H}\left( t\right)
\right\vert >l\right] }{\Pr \left[ \left\vert B\left( t\right) \right\vert >l%
\right] }= \\ 
\text{ }%
\end{array}%
$ & $0$ & $\infty $ \\ \hline
2) FIM vs. BM & $%
\begin{array}{c}
\text{ } \\ 
\lim_{l\rightarrow \infty }\frac{\Pr \left[ \left\vert I_{H}\left( t\right)
\right\vert >l\right] }{\Pr \left[ \left\vert B\left( t\right) \right\vert >l%
\right] }= \\ 
\text{ }%
\end{array}%
$ & $0$ & $\infty $ \\ \hline
3) FIM vs. FBM & $%
\begin{array}{c}
\text{ } \\ 
\lim_{l\rightarrow \infty }\frac{\Pr \left[ \left\vert I_{H}\left( t\right)
\right\vert >l\right] }{\Pr \left[ \left\vert B_{H}\left( t\right)
\right\vert >l\right] }= \\ 
\text{ }%
\end{array}%
$ & $0$ & $\infty $ \\ \hline
\end{tabular}

\ \ \ \ \ \ \ \ 

\textbf{Table 1}: tail-behavior comparisons. The comparisons are specified
with regard to the two anomalous-diffusion categories: sub-diffusion, ($0<H<%
\frac{1}{2}$); and super-diffusion ($\frac{1}{2}<H<1$). For sub-diffusion:
the tail behavior of FIM is infinitely `lighter' than the tail behavior of
FBM -- which, in turn, is infinitely `lighter' than the tail behavior of BM.
Conversely, for super-diffusion: the tail behavior of FIM is infinitely
`heavier' than the tail behavior of FBM -- which, in turn, is infinitely
`heavier' than the tail behavior of BM.

\subsection{\label{C}Statistical-heterogeneity comparison}

Inequality indices are quantitative gauges that -- using a
socioeconomic-inequality perspective -- score the statistical heterogeneity
of non-negative random variables with positive means \cite{StaEv}-\cite%
{TouIn}. We shall now use a particular inequality index in order to compare
the inherent statistical heterogeneity of FBM to the inherent statistical
heterogeneity of FIM.

Inequality indices are widely applied in economics and in the social
sciences to measure wealth inequalities in human societies \cite{Cou}-\cite%
{Cow}. An inequality index takes values in the unit interval, and with
regard to a human society under consideration: the lower the
inequality-index score -- the more egalitarian the distribution of wealth
among the society members; and the higher the inequality-index score -- the
less egalitarian the distribution of wealth among the society members. In
particular, the inequality index yields a zero score only when the society
under consideration is communist -- in which case all the society members
share a common (positive) wealth value. Inequality indices are invariant
with respect to the specific currency via which wealth is measured (e.g.
Dollar or Euro).

Given an inequality index $\mathcal{I}$, and given a non-negative random
variable $W$ with a positive mean, the statistical heterogeneity of the
random variable can be measured by the inequality index as follows \cite%
{TouIn}: deem the random variable $W$ to manifest the personal wealth value
of a member that is sampled at random from a virtual human society; then,
set the statistical-heterogeneity score of the random variable $W$ to be the
\ inequality-index score of the virtual human society. In what follows we
denote by $\mathcal{I}\left( W\right) $ the statistical-heterogeneity score
that the inequality index $\mathcal{I}$ assigns to the random variable $W$.

As $\mathcal{I}$ is an inequality index, the score $\mathcal{I}\left(
W\right) $ exhibits the two following properties. \textbf{I}) It takes
values in the unit interval, $0\leq \mathcal{I}\left( W\right) \leq 1$, and
it vanishes if and only if the random variable is deterministic: $\mathcal{I}%
\left( W\right) =0\Leftrightarrow W=const$ (with probability one). \textbf{II%
}) It is invariant with respect to changes of scale of the random variable: $%
\mathcal{I}\left( s\cdot W\right) =\mathcal{I}\left( W\right) $, where $s$
is any positive scale.

The particular inequality index $\mathcal{I}$ that we shall use is based on
the first-order and second-order moments of the random variable $W$, and its
statistical-heterogeneity score is \cite{TouIn},\cite{RenEq}: $\mathcal{I}%
\left( W\right) =1-\mathbf{E}[W]^{2}/\mathbf{E}[W^{2}]$. This
statistical-heterogeneity score is unique in the following sense: it is the
only score that has a one-to-one correspondence with the coefficient of
variation (CV) of the random variable $W$, i.e. the ratio of the
random-variable's standard deviation to the random-variable's mean. In
short, we henceforth term this statistical-heterogeneity score
\textquotedblleft CV score\textquotedblright .

Consider a selfsimilar one-dimensional random motion whose trajectory is $%
X\left( t\right) $ ($t\geq 0$), and whose Hurst exponent is $H$. We set the
fucus on the random variable $\left\vert X\left( t\right) \right\vert $ --
the motion's displacement, at time $t$, relative to its initial position.
The motion's selfsimilarity implies that the random variable $\left\vert
X\left( t\right) \right\vert $ is equal in law to the random variable $%
t^{H}\cdot \left\vert X\left( 1\right) \right\vert $. In turn, assuming that
the random variable $\left\vert X\left( 1\right) \right\vert $ has a
positive mean, the scaling property of inequality indices implies that: $%
\mathcal{I}\left( \left\vert X\left( t\right) \right\vert \right) =\mathcal{I%
}\left( \left\vert X\left( 1\right) \right\vert \right) $. Hence, the
quantity $\mathcal{I}\left( \left\vert X\left( 1\right) \right\vert \right) $
is, in effect, a statistical-heterogeneity score of the selfsimilar random
motion under consideration.

In the case of FBM the random variable $B_{H}\left( 1\right) $ is a Normal
random variable with zero mean. In turn, a probabilistic calculation that is
detailed in the Methods asserts that the CV\ score of FBM\ is%
\begin{equation}
\mathcal{I}\left[ \left\vert B_{H}\left( 1\right) \right\vert \right] =1-%
\frac{2}{\pi }\text{ .}  \label{41}
\end{equation}%
This CV\ score does not depend on the value of the Hurst exponent $H$.

In the case of FIM the statistical distribution of the random variable $%
I_{H}\left( 1\right) $ is governed by the density function of Eq. (\ref{32}%
). In turn, a probabilistic calculation that is detailed in the Methods
asserts that the CV\ score of FIM\ is%
\begin{equation}
\mathcal{I}\left[ \left\vert I_{H}\left( 1\right) \right\vert \right] =1-%
\frac{\sin \left( \pi H\right) }{\pi H}\text{ .}  \label{42}
\end{equation}%
This CV\ score depends on the value of the Hurst exponent $H$, and as a
function of the Hurst exponent $0<H<1$ (see Figure 7): it increases from the
zero lower-bound score $\lim_{H\rightarrow 0}\mathcal{I}\left[ \left\vert
I_{H}\left( 1\right) \right\vert \right] =0$ to the unit upper-bound score $%
\lim_{H\rightarrow 1}\mathcal{I}\left[ \left\vert I_{H}\left( 1\right)
\right\vert \right] =1$.

Evidently, the difference between the CV\ score of Eq. (\ref{41}) and the
CV\ score of Eq. (\ref{42}) is dramatic. On the one hand, changing the Hurst
exponent $H$ in the FBM model has no effect on the inherent statistical
heterogeneity. On the other hand, changing the Hurst exponent $H$ in the FIM
model has a profound effect on the inherent statistical heterogeneity --
which spans the full unit-interval range of statistical-heterogeneity scores.

\begin{figure}
\begin{centering}
\includegraphics[width=0.4\paperwidth,height=0.4\paperwidth,keepaspectratio]{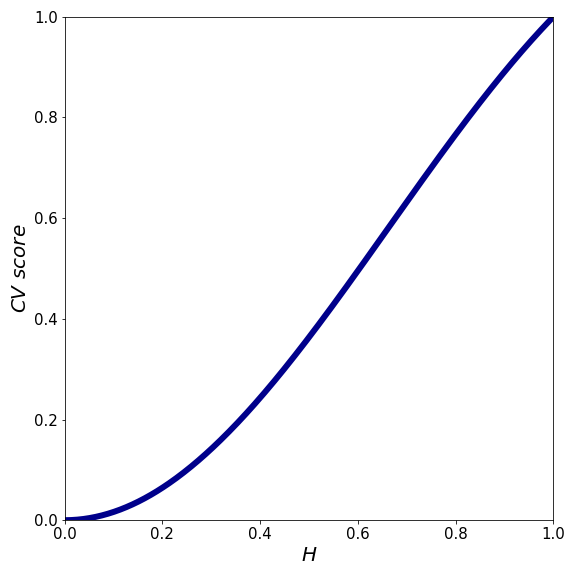}
\par\end{centering}
\caption{the FIM\ CV\ score of Eq. (\ref{42}), as a function of
the Hurst exponent $0<H<1$. Varying the value of the Hurst exponent, the
full unit-interval range of the CV score is attained.}

\end{figure}

\subsection{\label{D}Divergence comparison}

Consider a selfsimilar one-dimensional random motion whose trajectory is $%
X\left( t\right) $ ($t\geq 0$), and whose Hurst exponent is $H$. As noted
above, the motion's selfsimilarity implies that the random variable $X\left(
t\right) $ is equal in law to the random variable $t^{H}\cdot X\left(
1\right) $. Hence, the random variable $X\left( 1\right) $ -- the motion's
position at the time point $1$ -- is a `benchmark' for the motion's
positions at all times.

In the previous subsection we focused on quantifying the inherent
statistical heterogeneity of the random variable $\left\vert X\left(
t\right) \right\vert $ -- the motion's displacement, at time $t$, relative
to its initial position. In this subsection we shall focus on quantifying
the statistical divergence of the random variable $X\left( t\right) $ from
its benchmark -- the random variable $X\left( 1\right) $. To that end we use
the Kulback-Leibler (KL) divergence \cite{KL}-\cite{Kul}, which is arguably
the most widely applied measure of statistical divergence. Specifically, the
KL divergence of the random variable $X\left( t\right) $ from its benchmark $%
X\left( 1\right) $ is: $\mathcal{D}\left[ X\left( t\right) ||X\left(
1\right) \right] =\int_{-\infty }^{\infty }\ln \left[ f_{t}\left( x\right)
/f_{1}\left( x\right) \right] f_{t}\left( x\right) dx$, where $f_{t}\left(
x\right) $ is the density of the random variable $X\left( t\right) $ (and,
in particular, $f_{1}\left( x\right) $ is the density of the benchmark $%
X\left( 1\right) $).

With regard to FBM, a probabilistic calculation that is detailed in the
Methods asserts that the KL divergence of the position $B_{H}\left( t\right) 
$ from the benchmark position $B_{H}\left( 1\right) $ is 
\begin{equation}
\mathcal{D}\left[ B_{H}\left( t\right) ||B_{H}\left( 1\right) \right] =\frac{%
1}{2}\left( t^{2H}-1\right) -H\ln \left( t\right) \text{ .}  \label{71}
\end{equation}%
The asymptotic behavior of this KL divergence, in the temporal limit $%
t\rightarrow \infty $, is $\mathcal{D}\left[ B_{H}\left( t\right)
||B_{H}\left( 1\right) \right] \approx \frac{1}{2}t^{2H}$.

With regard to FIM, a probabilistic calculation that is detailed in the
Methods asserts that the KL divergence of the position $I_{H}\left( t\right) 
$ from the benchmark position $I_{H}\left( 1\right) $ is 
\begin{equation}
\mathcal{D}\left[ I_{H}\left( t\right) ||I_{H}\left( 1\right) \right]
=\left( 1-H\right) \left[ t-1-\ln \left( t\right) \right] \text{ .}
\label{72}
\end{equation}%
The asymptotic behavior of this KL divergence, in the temporal limit $%
t\rightarrow \infty $, is $\mathcal{D}\left[ I_{H}\left( t\right)
||I_{H}\left( 1\right) \right] \approx \left( 1-H\right) t$.

As a function of the temporal variable $t$, the KL divergences of Eqs. (\ref%
{71}) and (\ref{72}) share, qualitatively, a common U shape with: a unique
global minimum -- whose value is zero -- that is attained at the time point $%
1$. However, quantitatively, these KL divergences yield dramatically
different asymptotic behaviors. Using these asymptotic behaviors, we conduct
three divergence comparisons: FBM vs. BM; FIM vs. BM; and FIM vs. FBM. These
divergence comparisons are presented in Table 2, and they yield the
following `asymptotic orderings'. In the case of sub-diffusion ($0<H<\frac{1%
}{2}$): the FIM\ asymptotics are equivalent to the BM asymptotics -- which,
in turn, are infinitely `lighter' than the FBM asymptotics. In the case of
super-diffusion ($\frac{1}{2}<H<1$): the FIM\ asymptotics are equivalent to
the BM asymptotics -- which, in turn, are infinitely
`heavier\textquotedblright\ than the FBM asymptotics.

\ \ \ \ 

\begin{center}
{\large Table 2}
\end{center}

\ \ \ \ \ \ 

\begin{tabular}{|c|c|c|c|}
\hline
\textbf{Comparison} & \textbf{Limit} & $%
\begin{array}{c}
\text{\textbf{Sub diffusion}} \\ 
\text{(}0<H<\frac{1}{2}\text{)}%
\end{array}%
$ & $%
\begin{array}{c}
\text{\textbf{Super diffusion}} \\ 
\text{(}\frac{1}{2}<H<1\text{)}%
\end{array}%
$ \\ \hline
1) FBM vs. BM & $%
\begin{array}{c}
\text{ } \\ 
\lim_{t\rightarrow \infty }\frac{\mathcal{D}\left[ B_{H}\left( t\right)
||B_{H}\left( 1\right) \right] }{\mathcal{D}\left[ B\left( t\right)
||B\left( 1\right) \right] }= \\ 
\text{ }%
\end{array}%
$ & $0$ & $\infty $ \\ \hline
2) FIM vs. BM & $%
\begin{array}{c}
\text{ } \\ 
\lim_{t\rightarrow \infty }\frac{\mathcal{D}\left[ I_{H}\left( t\right)
||I_{H}\left( 1\right) \right] }{\mathcal{D}\left[ B\left( t\right)
||B\left( 1\right) \right] }= \\ 
\text{ }%
\end{array}%
$ & \multicolumn{2}{|c|}{$2\left( 1-H\right) $} \\ \hline
3) FIM vs. FBM & $%
\begin{array}{c}
\text{ } \\ 
\lim_{t\rightarrow \infty }\frac{\mathcal{D}\left[ I_{H}\left( t\right)
||I_{H}\left( 1\right) \right] }{\mathcal{D}\left[ B_{H}\left( t\right)
||B_{H}\left( 1\right) \right] }= \\ 
\text{ }%
\end{array}%
$ & $0$ & $\infty $ \\ \hline
\end{tabular}

\ \ \ \ \ 

\textbf{Table 2}: asymptotic KL-divergence comparisons. The comparisons are
specified with regard to the two anomalous-diffusion categories:
sub-diffusion ($0<H<\frac{1}{2}$); and super-diffusion ($\frac{1}{2}<H<1$).
For sub-diffusion: the FIM\ asymptotics -- which are equivalent to the BM
asymptotics -- are infinitely `lighter' than the FBM asymptotics.
Conversely, for super-diffusion: the FIM\ asymptotics -- which are
equivalent to the BM asymptotics -- are infinitely `heavier'\ than the FBM
asymptotics.

\subsection{\label{E}Aging comparison}

Consider a symmetric one-dimensional random motion whose trajectory is $%
X\left( t\right) $ ($t\geq 0$). The motion's displacement over the temporal
interval $[t,t+\Delta ]$ -- where $\Delta $ is the interval's positive
length -- is $\left\vert X\left( t+\Delta \right) -X\left( t\right)
\right\vert $. In turn, the motion's mean square displacement (MSD) over the
temporal interval $[t,t+\Delta ]$ is $\mathbf{E}[\left\vert X\left( t+\Delta
\right) -X\left( t\right) \right\vert ^{2}]$. And, as the motion is
symmetric, this MSD is equal to $\mathbf{Var}\left[ X\left( t+\Delta \right)
-X\left( t\right) \right] $ -- the variance of the increment $X\left(
t+\Delta \right) -X\left( t\right) $. In this subsection we shall compare
the variance of the FBM increment $B_{H}\left( t+\Delta \right) -B_{H}\left(
t\right) $ to the variance of the FIM increment $I_{H}\left( t+\Delta
\right) -I_{H}\left( t\right) $.

As FBM initiates at the spatial origin, and as the increments of FBM\ are
stationary, the FBM increment $B_{H}\left( t+\Delta \right) -B_{H}\left(
t\right) $ is equal in law to the random variable $B_{H}\left( \Delta
\right) $. In turn, as FBM is a selfsimilar process with Hurst exponent $H$,
the random variable $B_{H}\left( \Delta \right) $ is equal in law to the
random variable $\Delta ^{H}\cdot B_{H}\left( 1\right) $. Consequently, the
FBM increment $B_{H}\left( t+\Delta \right) -B_{H}\left( t\right) $ is a
Normal random variable with mean zero and with variance%
\begin{equation}
\mathbf{Var}\left[ B_{H}\left( t+\Delta \right) -B_{H}\left( t\right) \right]
=\mathbf{Var}\left[ B_{H}\left( 1\right) \right] \cdot \Delta ^{2H}\text{ .}
\label{51}
\end{equation}%
Due to the fact that the increments of FBM\ are stationary, the variance of
Eq. (\ref{51}) depends only on the length $\Delta $ of the temporal interval 
$[t,t+\Delta ]$.

Using the fact that FIM\ is a selfsimilar process with Hurst exponent $H$,
as well as the fact that FIM is a martingale, it is shown in the Methods
that the FIM increment $I_{H}\left( t+\Delta \right) -I_{H}\left( t\right) $
is a random variable with mean zero and with variance 
\begin{equation}
\mathbf{Var}\left[ I_{H}\left( t+\Delta \right) -I_{H}\left( t\right) \right]
=\mathbf{Var}\left[ I_{H}\left( 1\right) \right] \cdot \left[ \left(
t+\Delta \right) ^{2H}-t^{2H}\right] \text{ .}  \label{52}
\end{equation}%
The variance of Eq. (\ref{52}) depends on the starting point $t$, as well as
on the length $\Delta $, of the temporal interval $[t,t+\Delta ]$. Hence,
Eq. (\ref{52}) implies that the increments of FIM\ are \emph{not} stationary.

The asymptotic behavior of the variance of Eq. (\ref{52}), in the temporal
limit $t\rightarrow \infty $, is $\mathbf{Var}\left[ I_{H}\left( t+\Delta
\right) -I_{H}\left( t\right) \right] \approx v\cdot t^{2H-1}$, where $v=%
\mathbf{Var}\left[ I_{H}\left( 1\right) \right] \cdot 2\Delta H$ (see the
Methods for the details of this asymptotic equivalence). Consequently, for
sub-diffusion ($0<H<\frac{1}{2}$) Eq. (\ref{52}) implies that: $\lim_{t\rightarrow \infty }\mathbf{Var}%
\left[ I_{H}\left( t+\Delta \right) -I_{H}\left( t\right) \right] =0$;
namely, the variance of the FIM increment $I_{H}\left( t+\Delta \right)
-I_{H}\left( t\right) $ vanishes in the temporal limit $t\rightarrow \infty $. And, for super-diffusion ($\frac{1}{2}<H<1$) Eq. (\ref{52}) implies that: $\lim_{t\rightarrow \infty }%
\mathbf{Var}\left[ I_{H}\left( t+\Delta \right) -I_{H}\left( t\right) \right]
=\infty $; namely, the variance of the FIM $I_{H}\left( t+\Delta \right)
-I_{H}\left( t\right) $ increment explodes in the temporal limit $%
t\rightarrow \infty $.

Evidently, the difference between the variance of Eq. (\ref{51}) and the
variance of Eq. (\ref{52}) is dramatic. Indeed, consider an observation that
starts at the time point $t$, and whose length is $\Delta $; also, address $%
t $ as the `age' of the observation. On the one hand, the `age' of the
observation has no effect when observing BM and FBM. On the other hand, the
`age' of the observation has a profound effect when observing FIM: the
statistical fluctuations of the FIM increment become smaller and smaller
with `age' when the FIM is sub-diffusive; and the statistical fluctuations
of the FIM increment become larger and larger with `age' when the FIM is
super-diffusive.

\subsection{\label{F}Correlation comparison}

In this last comparison between FBM and FIM, we examine the correlations of
their velocities. To that end we set two positive and distinct time points, $%
t_{1}$ and $t_{2}$, and address the covariance of their velocities at these
points.

As FBM is a selfsimilar process with Hurst exponent $H$, and as its
increments are stationary, it follows that the covariance of the FBM\
positions $B_{H}\left( t_{1}\right) $ and $B_{H}\left( t_{2}\right) $ is 
\cite{EM}: 
\begin{equation}
\mathbf{Cov}\left[ B_{H}\left( t_{1}\right) ,B_{H}\left( t_{2}\right) \right]
=\frac{1}{2}\mathbf{Var}\left[ B_{H}\left( 1\right) \right] \cdot \left(
t_{1}^{2H}-\left\vert t_{1}-t_{2}\right\vert ^{2H}+t_{2}^{2H}\right) \text{ .%
}  \label{61}
\end{equation}%
Differentiating the covariance of Eq. (\ref{61}) with respect to the
temporal variable $t_{1}$, and then with respect to the temporal variable $%
t_{2}$, implies that the covariance of the FBM\ velocities $\dot{B}%
_{H}\left( t_{1}\right) $ and $\dot{B}_{H}\left( t_{2}\right) $ is:%
\begin{equation}
\mathbf{Cov}\left[ \dot{B}_{H}\left( t_{1}\right) ,\dot{B}_{H}\left(
t_{2}\right) \right] =\mathbf{Var}\left[ B_{H}\left( 1\right) \right] \cdot 
\frac{H\left( 2H-1\right) }{\left\vert t_{1}-t_{2}\right\vert ^{2\left(
1-H\right) }}\text{ .}  \label{62}
\end{equation}

The covariance of Eq. (\ref{62}) has the following implications. For
sub-diffusion ($0<H<\frac{1}{2}$) the FBM velocities are negatively correlated: $\mathbf{Cov}[\dot{B}%
_{H}\left( t_{1}\right) ,\dot{B}_{H}\left( t_{2}\right) ]<0$. For super-diffusion ($\frac{1}{2}<H<1$) the FBM velocities are positively correlated: $\mathbf{%
Cov}[\dot{B}_{H}\left( t_{1}\right) ,\dot{B}_{H}\left( t_{2}\right) ]>0$. And, for regular diffusion ($H=\frac{1}{2}$) -- in which case FBM is BM -- the FBM velocities are uncorrelated: $\mathbf{Cov}[\dot{B}_{H}\left( t_{1}\right) ,\dot{B}_{H}\left( t_{2}\right)]=0$. 

As FIM is a selfsimilar process with Hurst exponent $H$, and as it is a
martingale, it is shown in the Methods (see Eq. (\ref{M35})) that the
covariance of the FIM\ positions $I_{H}\left( t_{1}\right) $ and $%
I_{H}\left( t_{2}\right) $ is:%
\begin{equation}
\mathbf{Cov}\left[ I_{H}\left( t_{1}\right) ,I_{H}\left( t_{2}\right) \right]
=\mathbf{Var}\left[ I_{H}\left( 1\right) \right] \cdot \left( \min \left\{
t_{1},t_{2}\right\} \right) ^{2H}\text{ .}  \label{63}
\end{equation}%
Differentiating the covariance of Eq. (\ref{63}) with respect to the
temporal variable $t_{1}$, and then with respect to the temporal variable $%
t_{2}$, implies that the covariance of the FIM\ velocities $\dot{I}%
_{H}\left( t_{1}\right) $ and $\dot{I}_{H}\left( t_{2}\right) $ is zero:%
\begin{equation}
\mathbf{Cov}\left[ \dot{I}_{H}\left( t_{1}\right) ,\dot{I}_{H}\left(
t_{2}\right) \right] =0\text{ .}  \label{64}
\end{equation}%
Namely, the FIM\ velocities $\dot{I}_{H}\left( t_{1}\right) $ and $\dot{I}%
_{H}\left( t_{2}\right) $ are uncorrelated.

Evidently, the difference between the velocities covariance of Eq. (\ref{62}%
) and the velocities covariance of Eq. (\ref{64}) is dramatic. On the one
hand, changing the Hurst exponent $H$ in the FBM model has a profound
effect, both qualitatively and quantitatively, on the correlations of the
FBM\ velocities. On the other hand, changing the Hurst exponent $H$ in the
FIM model has no effect on the correlations of the FIM\ velocities -- as
these velocities are always uncorrelated.

\section{\label{3}Diffusion in a logarithmic potential}

As described above, the Ito stochastic dynamics of FIM have a zero drift and
a position-dependent volatility. Namely, the evolution of FIM is governed by
an Ito stochastic differential equation (SDE) of the form $\dot{X}\left(
t\right) =\sigma \lbrack X\left( t\right) ]\dot{B}\left( t\right) $. Many
researchers in the physical sciences and in engineering are used to work
with Langevin stochastic dynamics \cite{Lan}-\cite{Pav}, which have a
position-dependent drift and a constant volatility. Up to an inherent time
scale, the latter evolution is governed by a Langevin SDE of the form $\dot{X%
}\left( t\right) =\mu \left[ X\left( t\right) \right] +\dot{B}\left(
t\right) $. 

As noted above, a special case of the FIM model was investigated in \cite{LB}: in this special case the motion runs over the positive half-line ($0<x<\infty $), and the Hurst exponent is in the range $\frac{1}{2} \leq H <1$. It was shown in \cite{LB} that the special case of FIM admits a Langevin representation. In this section we shall present a Langevin representation for the (general) FIM model. 

There is a one-to-one mapping between FIM and a particular random motion
that is generated by specific Langevin stochastic dynamics. Indeed, map the
FIM trajectory $I_{H}\left( t\right) $ ($t\geq 0$) to the trajectory $\xi
_{H}\left( t\right) $ ($t\geq 0$) via the following transformation: $\xi
_{H}\left( t\right) =\varphi \lbrack I_{H}\left( t\right) ]$, where $\varphi
\left( x\right) =2H\left\vert x\right\vert ^{1/2H}sign\left( x\right) $, and
where $sign\left( x\right) $ is the sign of $x$.\footnote{%
Namely, $sign\left( x\right) =1$ if $x$ is positive, and $sign\left(
x\right) =-1$ if $x$ is negative.} Then, using Ito's formula and Eq. (\ref%
{21}), it is shown in the Methods that the stochastic dynamics of the random
motion $\xi _{H}\left( t\right) $ are governed by the Langevin SDE 
\begin{equation}
\dot{\xi}_{H}\left( t\right) =\left( \frac{1}{2}-H\right) \frac{1}{\xi
_{H}\left( t\right) }+\dot{B}\left( t\right) \text{ .}  \label{80}
\end{equation}%
Shifting back from the trajectory $\xi _{H}\left( t\right) $ ($t\geq 0$) to
the FIM trajectory $I_{H}\left( t\right) $ ($t\geq 0$) is via the
corresponding inverse transformation: $I_{H}\left( t\right) =\varphi ^{-1}%
\left[ \xi _{H}\left( t\right) \right] $, where $\varphi ^{-1}\left(
y\right) =(\frac{1}{2H}\left\vert y\right\vert )^{2H}sign\left( y\right) $.

The drift of the Langevin SDE (\ref{80}) is harmonic, $\mu \left( x\right) =(%
\frac{1}{2}-H)/x$. This drift underscores and highlights the three different
diffusivity categories of FIM -- which are determined by the value of its
Hurst exponent $H$. When FIM is sub-diffusive ($0<H<\frac{1}{2}$) then the Langevin drift\ is
`repulsive': it pushes \emph{away} from the spatial origin. When FIM is super-diffusive ($\frac{1}{2}<H<1$) then the Langevin drift is `attractive': it pushes \emph{towards}
the spatial origin. And when FIM is a regular diffusion ($H=\frac{1}{2}$) then the Langevin drift is zero.

A Langevin SDE of the form $\dot{X}\left( t\right) =\mu \left[ X\left(
t\right) \right] +\dot{B}\left( t\right) $ is commonly characterized by its
\textquotedblleft potential function\textquotedblright : a function $U\left(
x\right) $ whose negative gradient is the drift, $-U^{\prime }\left(
x\right) =\mu \left( x\right) $. The potential function $U\left( x\right) $
manifests a `potential landscape' that underpins and governs the Langevin
stochastic dynamics. The random motions that are generated by a Langevin SDE
with a logarithmic potential $U\left( x\right) =r\cdot \ln \left( \left\vert
x\right\vert \right) $ -- where $r$ is a real parameter -- are termed
\textquotedblleft diffusion in a logarithmic potential\textquotedblright\
(DLP). The topic of DLP attracted significant scientific interest recently 
\cite{KB}-\cite{PM}.

Arguably, the best known example of DLP is the Bessel process \cite{Pit}-%
\cite{Kat}: the Euclidean distance of a Brownian motion -- which takes place
in the $d$-dimensional Euclidean space -- from the space's origin. In this
example DLP runs over the positive half-line ($0<x<\infty $), and its
parameter is $r=\frac{1}{2}\left( d-1\right) $. So, for dimensions $%
d=2,3,4,5,6,\cdots $, the Bessel process yields the DLP parameter values $r=$
$\frac{1}{2},\frac{2}{2},\frac{3}{2},\frac{4}{2},\frac{5}{2},\cdots $.

FIM\ yields -- via the Langevin SDE (\ref{80}) -- another example of DLP.
Indeed, the underpinning potential function of the Langevin SDE (\ref{80})
is logarithmic with parameter $r=H-\frac{1}{2}$. As the Hurst exponent of
FIM takes values in the range $0<H<1$, the corresponding DLP parameter takes
values in the range $-\frac{1}{2}<r<\frac{1}{2}$. Negative values of the DLP
parameter correspond to sub-diffusive FIM, and positive values of the DLP
parameter correspond to super-diffusive FIM.

As noted above, the potential function $U\left( x\right) $ manifests the
`potential landscape' that underpins and governs the Langevin SDE $\dot{X}%
\left( t\right) =\mu \left[ X\left( t\right) \right] +\dot{B}\left( t\right) 
$. Analogously, the volatility function $\sigma \left( x\right) $ manifests
the `volatility landscape' that underpins and governs the drift-less Ito SDE 
$\dot{X}\left( t\right) =\sigma \lbrack X\left( t\right) ]\dot{B}\left(
t\right) $. The one-to-one mapping between FIM and its corresponding DLP
induces a one-to-one mapping between their respective underpinning
landscapes: the power-law volatility landscape $\sigma \left( x\right)
=\left\vert x\right\vert ^{1-\frac{1}{2H}}$, and the logarithmic potential
landscape $U\left( x\right) =(H-\frac{1}{2})\ln \left( \left\vert
x\right\vert \right) $. Schematic illustrations of these corresponding
landscapes are depicted in Figure 8.

\ \

\begin{figure}
\begin{centering}
\includegraphics[width=0.4\paperwidth,height=0.3\paperwidth,keepaspectratio]{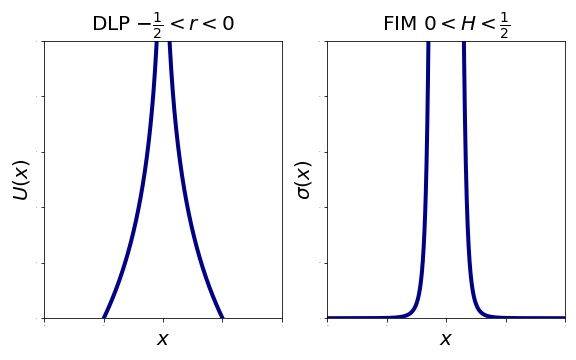}
\includegraphics[width=0.4\paperwidth,height=0.3\paperwidth,keepaspectratio]{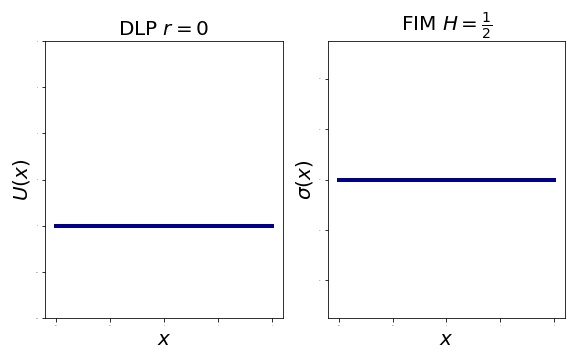}
\includegraphics[width=0.4\paperwidth,height=0.3\paperwidth,keepaspectratio]{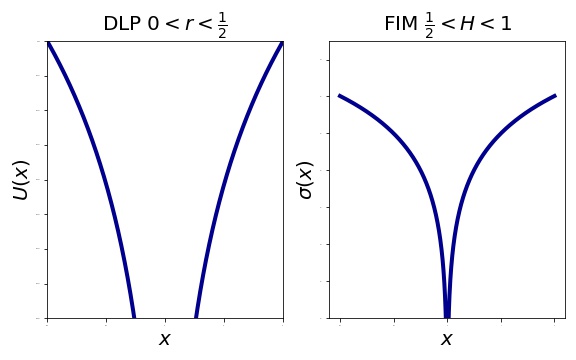}
\par\end{centering}
\caption{schematic illustrations of the shapes of the FIM volatility landscape, and of the corresponding DLP potential landscape. Depicted in the right panels, the FIM volatility landscape is: $\sigma \left( x\right) =\left\vert x\right\vert ^{1-\frac{1}{2H%
}}$ over the spatial axis $-\infty <x<\infty $. Depicted in the left panels, the DLP potential landscape is: $U\left( x\right) =r\cdot \ln \left( \left\vert
x\right\vert \right)$ over the spatial axis $-\infty <x<\infty $. The one-to-one correspondence between the two landscapes is via the relation $r=H-\frac{1}{2}$, where $0<H<1$ and $-\frac{1}{2}<r<\frac{1}{2}$. The two landscapes display markedly different shapes in the
three different diffusivity categories: sub-diffusion $0<H<\frac{1}{2}$ and $-\frac{1}{2}<r<0$ (top panels); super-diffusion $\frac{1}{2}<H<1$ and $0<r<\frac{1}{2}$ (bottom panels); and regular diffusion $H=\frac{1}{2}$ and $r=0$ (middle panels).}

\end{figure}

\section{\label{4}Summary}

Setting off from the paradigmatic model for regular diffusion, \emph{%
Brownian motion} (BM), in this paper we examined and compared two
anomalous-diffusion models that generalize BM: the popular and well-known 
\emph{fractional Brownian motion} (FBM) model, and the \emph{%
fractional Ito motion} (FIM) model. A `bird's-eye view' of FBM vs. FIM is
summarized in Table 3. More specifically, going over a list of key BM\
properties, Table 3: compares FBM\ and FIM to BM; and compares FBM\ to FIM.

From the perspective of mean-square-displacement measurements, both FBM and
FIM exhibit the very same modeling `efficacy'. Indeed, in both these models
the range of the Hurst exponent is $0<H<1$. Consequently, both these models
yield: sub-diffusion in the Hurst range $0<H<\frac{1}{2}$; super-diffusion
in the Hurst range $\frac{1}{2}<H<1$; and regular diffusion at the Hurst
value $H=\frac{1}{2}$ -- in which case both FBM\ and FIM are BM. Thus, the
FBM\ and FIM models are indistinguishable from the mean-square-displacement
perspective. 

However, the stochastic propagation mechanisms of FBM and FIM are
dramatically different. On the one hand, the FBM\ velocity (Eq. (\ref{11}))
is a linear moving average of the underlying white noise -- \ where the
averaging is temporal, and it goes from the present (the time point $t$) all
the way back to the infinitely distant past (the time point $-\infty $). We
emphasize that the FBM\ velocity does not depend on the FBM positions. On
the other hand, the FIM velocity (Eq. (\ref{21})) is a non-linear mapping of
the following present (the time point $t$) quantities: the present FIM
position, and the present value of the underlying white noise.

The profound differences between FBM and FIM -- which are evident from Table
3, and moreover from the detailed comparisons of section \ref{2} -- stem
from the dramatic differences between the stochastic propagation mechanisms
of these anomalous-diffusion models. Due to its properties, FBM\ is a rather
challenging model from the perspectives of numerical simulation, of
stochastic integration, and of analytic tractability. Conversely, from these
very perspectives, FIM\ is a compelling and an amicable working model.

In addition to FIM\ being a potent model for anomalous diffusion, we showed
that there is a one-to-one mapping between FIM and diffusion in a
logarithmic potential (where the coefficient of the logarithmic potential is 
$H-\frac{1}{2}$). With its compelling properties, its easy numerical simulation, its analytic tractability, and its linkage to diffusion in
a logarithmic potential -- we believe that FIM may become a prevalent
working model for anomalous diffusion in various fields of science and
engineering.

\ \ 

\begin{center}
{\large Table 3}
\end{center}

\ \ \ \ \ 

\begin{tabular}{|c|c|c|c|}
\hline
& \textbf{BM} & \textbf{FBM} & \textbf{FIM} \\ \hline
Finite variance & \multicolumn{3}{|c|}{YES} \\ \hline
Symmetric process & \multicolumn{3}{|c|}{YES} \\ \hline
Continuous trajectory & \multicolumn{3}{|c|}{YES} \\ \hline
Selfsimilar process & \multicolumn{3}{|c|}{YES} \\ \hline
Hurst exponent & $H=\frac{1}{2}$ & $0<H<1$ & $0<H<1$ \\ \hline
Gaussian process & YES & YES & NO \\ \hline
Markov process & YES & NO & YES \\ \hline
Martingale & YES & NO & YES \\ \hline
Stationary increments & YES & YES & NO \\ \hline
Uncorrelated increments & YES & NO & YES \\ \hline
\end{tabular}%
\ \ \ \ \ \ 

\ \ 

\textbf{Table 3}: a `bird's-eye view' comparison between the FBM and FIM
models with Hurst exponents $H\neq \frac{1}{2}$. The comparison is with
respect to the `benchmark' BM model -- which the FBM and FIM models yield at
the Hurst exponent $H=\frac{1}{2}$. Each row of the table specifies a
different key\ property of BM. The top four rows of the table highlight the
properties that FBM and FIM share in common with BM. The bottom six rows of
the table underscore the marked differences between the three models (the
last row of the table refers to non-overlapping increments).

\newpage

\section{Methods}

\subsection{The macroscopic power-law structure of the MSD function}

Consider a one-dimensional random motion whose trajectory is $X\left(
t\right) $ ($t\geq 0$). Given a positive number $n$, speed up the temporal
scale by the factor $n$. To compensate for speeding up time, shrink the
spatial scale by the positive factor $\varphi \left( n\right) $ (which is a
continuous function of the positive number $n$). This spatio-temporal scaling shifts
the trajectory $X\left( t\right) $ ($t\geq 0$) to the trajectory $%
X_{n}\left( t\right) $ ($t\geq 0$), where 
\begin{equation}
X_{n}\left( t\right) =\frac{1}{\varphi \left( n\right) }X\left( nt\right) 
\text{ .}  \label{M10}
\end{equation}

The MSD function of the random motion is $\phi \left( t\right) =\mathbf{E}%
[\left\vert X\left( t\right) -X\left( 0\right) \right\vert ^{2}]$. In turn,
the MSD function of the scaled random motion is%
\begin{equation}
\left. 
\begin{array}{l}
\phi _{n}\left( t\right) =\mathbf{E}\left[ \left\vert X_{n}\left( t\right)
-X_{n}\left( 0\right) \right\vert ^{2}\right] \\ 
\text{ } \\ 
=\mathbf{E}\left[ \left\vert \frac{1}{\varphi \left( n\right) }X\left(
nt\right) -\frac{1}{\varphi \left( n\right) }X\left( n0\right) \right\vert
^{2}\right] \\ 
\text{ } \\ 
=\frac{1}{\varphi \left( n\right) ^{2}}\mathbf{E}\left[ \left\vert X\left(
nt\right) -X\left( 0\right) \right\vert ^{2}\right] \\ 
\text{ } \\ 
=\frac{\phi \left( nt\right) }{\varphi \left( n\right) ^{2}}\text{ .}%
\end{array}%
\right.  \label{M11}
\end{equation}%
This MSD function can be written as follows 
\begin{equation}
\phi _{n}\left( t\right) =\frac{\phi \left( n\right) }{\varphi \left(
n\right) ^{2}}\cdot \frac{\phi \left( nt\right) }{\phi \left( n\right) }%
\text{ .}  \label{M12}
\end{equation}

Now, consider the MSD function $\phi _{n}\left( t\right) $ in the scaling
limit $n\rightarrow \infty $. In general, in order to obtain a non-trivial
point-wise limit $\phi _{\infty }\left( t\right) =\lim_{n\rightarrow \infty
}\phi _{n}\left( t\right) $ (for any positive time $t$), two conditions must
be met. Firstly, the spatial scaling function $\varphi \left( n\right) $ should be asymptotically equivalent (in the limit $n\rightarrow \infty $) to
the function $\sqrt{\phi \left( n\right) }$. Secondly, the MSD function $%
\phi \left( t\right) $ should be regularly varying at infinity \cite{BGT}. If
these two conditions hold then the non-trivial point-wise limit is 
\begin{equation}
\phi _{\infty }\left( t\right) =c\cdot t^{\epsilon }\text{ ,}  \label{M13}
\end{equation}%
where: $c$ is a positive constant that stems from the first condition; and $%
\epsilon $ is a positive regular-variation exponent that stems from the second condition.

\subsection{The CV scores of FBM and FIM}

We start with the random variable%
\begin{equation}
Z=2H^{2}\cdot \left\vert I_{H}\left( 1\right) \right\vert ^{1/H}\text{ .}
\label{M20}
\end{equation}%
The cumulative distribution function of the random variable $Z$ is 
\begin{equation}
\left. 
\begin{array}{l}
\Pr \left( Z\leq z\right) =\Pr \left( 2H^{2}\cdot \left\vert I_{H}\left(
1\right) \right\vert ^{1/H}\leq z\right)  \\ 
\text{ } \\ 
=\Pr \left[ \left\vert I_{H}\left( 1\right) \right\vert \leq \left( \frac{z}{%
2H^{2}}\right) ^{H}\right]  \\ 
\text{ } \\ 
=2\Pr \left[ I_{H}\left( 1\right) \leq \left( \frac{z}{2H^{2}}\right) ^{H}%
\right] -1%
\end{array}%
\right.   \label{M21}
\end{equation}%
($z\geq 0$). In the transition from the second line of Eq. (\ref{M21}) to
the third line we used the fact that the random variable $I_{H}\left(
1\right) $ is symmetric (i.e., the random variable $I_{H}\left( 1\right) $
is equal in law to the random variable $-I_{H}\left( 1\right) $).

Denote by $g\left( z\right) $ ($z\geq 0$) the density function of the random
variable $Z$, and denote by $f\left( x\right) $ ($-\infty <x<\infty $) the
density function of the random variable $I_{H}\left( 1\right) $.
Differentiating Eq. (\ref{M21}), and then using the density function of Eq. (%
\ref{32}), implies that%
\begin{equation}
\left. 
\begin{array}{l}
g\left( z\right) =2f\left[ \left( \frac{z}{2H^{2}}\right) ^{H}\right] \frac{1%
}{\left( 2H^{2}\right) ^{H}}Hz^{H-1} \\ 
\text{ } \\ 
=\frac{2H}{\left( 2H^{2}\right) ^{H}}f\left[ \left( \frac{z}{2H^{2}}\right)
^{H}\right] z^{H-1}\text{ } \\ 
\text{ } \\ 
=\frac{2H}{\left( 2H^{2}\right) ^{H}}\left\{ \frac{H^{1-2H}}{2^{H}\Gamma
\left( 1-H\right) }\exp \left[ -2H^{2}\left\vert \left( \frac{z}{2H^{2}}%
\right) ^{H}\right\vert ^{\frac{1}{H}}\right] \left\vert \left( \frac{z}{%
2H^{2}}\right) ^{H}\right\vert ^{\frac{1}{H}-2}\right\} z^{H-1} \\ 
\text{ } \\ 
=\left\{ \frac{2H}{\left( 2H^{2}\right) ^{H}}\frac{H^{1-2H}}{2^{H}\Gamma
\left( 1-H\right) }\frac{1}{\left( 2H^{2}\right) ^{1-2H}}\right\} \exp
\left( -z\right) \left\{ z^{1-2H}z^{H-1}\right\} \\ 
\text{ } \\ 
=\frac{1}{\Gamma \left( 1-H\right) }\exp \left( -z\right) z^{-H}\text{ .}%
\end{array}%
\right.  \label{M23}
\end{equation}%
In turn, Eq. (\ref{M23}) implies that 
\begin{equation}
g\left( z\right) =\frac{1}{\Gamma \left( 1-H\right) }\exp \left( -z\right)
z^{\left( 1-H\right) -1}\text{ .}  \label{M24}
\end{equation}%
The density function of Eq. (\ref{M24}) characterizes a Gamma distribution
with Gamma exponent $1-H$. Consequently, the mean of the random variable $%
Z^{p}$ -- where $p$ is a positive power -- is%
\begin{equation}
\mathbf{E}\left[ Z^{p}\right] =\frac{\Gamma \left( 1-H+p\right) }{\Gamma
\left( 1-H\right) }\text{ .}  \label{M25}
\end{equation}

Eq. (\ref{M20}) implies that 
\begin{equation}
\left\vert I_{H}\left( 1\right) \right\vert =\left( \frac{1}{2H^{2}}\right)
^{H}\cdot Z^{H}\text{ .}  \label{M26}
\end{equation}%
In turn, for any inequality index $\mathcal{I}$ we have 
\begin{equation}
\mathcal{I}\left( \left\vert I_{H}\left( 1\right) \right\vert \right) =%
\mathcal{I}\left( Z^{H}\right) \text{ .}  \label{M27}
\end{equation}%
In particular, for the inequality index $\mathcal{I}$ that yields the CV
score we have%
\begin{equation}
\left. 
\begin{array}{l}
\mathcal{I}\left( \left\vert I_{H}\left( 1\right) \right\vert \right) =1-%
\frac{\mathbf{E}[Z^{H}]^{2}}{\mathbf{E}[Z^{2H}]} \\ 
\text{ } \\ 
=1-\frac{\left[ \frac{\Gamma \left( 1-H+H\right) }{\Gamma \left( 1-H\right) }%
\right] ^{2}}{\left[ \frac{\Gamma \left( 1-H+2H\right) }{\Gamma \left(
1-H\right) }\right] }=1-\frac{\Gamma \left( 1\right) ^{2}}{\Gamma \left(
1+H\right) \Gamma \left( 1-H\right) } \\ 
\text{ } \\ 
=1-\frac{1}{H\Gamma \left( H\right) \cdot \Gamma \left( 1-H\right) }=1-\frac{%
1}{H}\frac{1}{\Gamma \left( H\right) \Gamma \left( 1-H\right) } \\ 
\text{ } \\ 
=1-\frac{1}{H}\frac{\sin \left( \pi H\right) }{\pi }=1-\frac{\sin \left( \pi
H\right) }{\pi H}\text{ .}%
\end{array}%
\right.   \label{M28}
\end{equation}%
In Eq. (\ref{M28}) we used the following (well known) properties of the
Gamma function: $\Gamma \left( 1\right) =1$; $\Gamma \left( 1+H\right)
=H\Gamma \left( H\right) $; and $\Gamma \left( H\right) \Gamma \left(
1-H\right) =\pi /\sin \left( \pi H\right) $. Eq. (\ref{M28}) proves Eq. (\ref%
{42}).

As noted above, at the Hurst exponent $H=\frac{1}{2}$ the random variable $%
I_{1/2}\left( 1\right) $ is \textquotedblleft standard
Normal\textquotedblright , i.e. it is a Normal random variable with mean
zero and with variance one. Also, as noted above, the random variable $%
B_{H}\left( 1\right) $ is Normal with mean zero and with variance $b=\mathbf{%
Var}\left[ B_{H}\left( 1\right) \right] $. Hence, the random variables $%
I_{1/2}\left( 1\right) $ and $B_{H}\left( 1\right) $ differ only by a scale
factor:%
\begin{equation}
B_{H}\left( 1\right) =\sqrt{b}\cdot I_{1/2}\left( 1\right) \text{ ,}
\label{M290}
\end{equation}%
where the equality is in law. In turn, any inequality index $\mathcal{I}$
will assign the random variables $\left\vert I_{1/2}\left( 1\right)
\right\vert $ and $\left\vert B_{H}\left( 1\right) \right\vert $ the same
statistical-heterogeneity score. In particular, for the CV score Eq. (\ref%
{M28}) implies that%
\begin{equation}
\left. 
\begin{array}{l}
\mathcal{I}\left( \left\vert B_{H}\left( 1\right) \right\vert \right) =%
\mathcal{I}\left( \left\vert I_{1/2}\left( 1\right) \right\vert \right)  \\ 
\text{ } \\ 
=1-\frac{\sin \left( \pi \frac{1}{2}\right) }{\pi \frac{1}{2}}=1-\frac{2}{%
\pi }\text{ .}%
\end{array}%
\right.   \label{M29}
\end{equation}%
Eq. (\ref{M29}) proves Eq. (\ref{41}).

\subsection{The Kulback-Leibler divergence of FBM and FIM}

Consider a one-dimensional random motion whose trajectory is $X\left(
t\right) $ ($t\geq 0$), and denote by $f_{t}\left( x\right) $ the density of
the random variable $X\left( t\right) $. As noted above, the Kulback-Leibler
(KL) divergence of the random variable $X\left( t\right) $ from the random
variable $X\left( 1\right) $ is: 
\begin{equation}
\mathcal{D}\left[ X\left( t\right) ||X\left( 1\right) \right] =\int_{-\infty
}^{\infty }\ln \left[ \frac{f_{t}\left( x\right) }{f_{1}\left( x\right) }%
\right] f_{t}\left( x\right) dx\text{ .}  \label{M40}
\end{equation}

\subsubsection{FBM case}

For FBM, Eq. (\ref{31}) implies that%
\begin{equation}
\left. 
\begin{array}{l}
\frac{f_{t}\left( x\right) }{f_{1}\left( x\right) }=\frac{\frac{1}{\sqrt{%
2\pi b}}\frac{1}{t^{H}}\exp \left( -\frac{x^{2}}{2bt^{2H}}\right) }{\frac{1}{%
\sqrt{2\pi b}}\exp \left( -\frac{x^{2}}{2b}\right) } \\ 
\text{ } \\ 
=\frac{1}{t^{H}}\exp \left[ \frac{x^{2}}{2b}\left( 1-\frac{1}{t^{2H}}\right) %
\right] \text{ ,}%
\end{array}%
\right.   \label{M41}
\end{equation}%
and hence%
\begin{equation}
\ln \left[ \frac{f_{t}\left( x\right) }{f_{1}\left( x\right) }\right] =-H\ln
\left( t\right) +\frac{1}{2b}\left( 1-\frac{1}{t^{2H}}\right) \cdot x^{2}%
\text{ .}  \label{M42}
\end{equation}%
Setting the random motion $X\left( t\right) $ ($t\geq 0$) to be FBM, and
substituting Eq. (\ref{M42}) into Eq. (\ref{M40}) yields%
\begin{equation}
\left. 
\begin{array}{l}
\mathcal{D}\left[ B_{H}\left( t\right) ||B_{H}\left( 1\right) \right]  \\ 
\text{ } \\ 
=\int_{-\infty }^{\infty }\left[ -H\ln \left( t\right) +\frac{1}{2b}\left( 1-%
\frac{1}{t^{2H}}\right) \cdot x^{2}\right] f_{t}\left( x\right) dx \\ 
\text{ } \\ 
=-H\ln \left( t\right) \int_{-\infty }^{\infty }f_{t}\left( x\right) dx+%
\frac{1}{2b}\left( 1-\frac{1}{t^{2H}}\right) \int_{-\infty }^{\infty
}x^{2}f_{t}\left( x\right) dx \\ 
\text{ } \\ 
=-H\ln \left( t\right) +\frac{1}{2b}\left( 1-\frac{1}{t^{2H}}\right) \mathbf{%
E}\left[ B_{H}\left( t\right) ^{2}\right] \text{ .}%
\end{array}%
\right.   \label{M43}
\end{equation}%
As noted above, the properties of FBM imply that the random variable $%
B_{H}\left( t\right) $ is Normal with mean zero with variance $\mathbf{Var}%
\left[ B_{H}\left( t\right) \right] =b\cdot t^{2H}$ -- and hence%
\begin{equation}
\mathbf{E}\left[ B_{H}\left( t\right) ^{2}\right] =b\cdot t^{2H}\text{ .}
\label{M44}
\end{equation}%
Combining Eq. (\ref{M43}) and Eq. (\ref{M44}) \ together, we obtain that

\begin{equation}
\left. 
\begin{array}{l}
\mathcal{D}\left[ B_{H}\left( t\right) ||B_{H}\left( 1\right) \right] \\ 
\text{ } \\ 
=-H\ln \left( t\right) +\frac{1}{2b}\left( 1-\frac{1}{t^{2H}}\right) bt^{2H}
\\ 
\text{ } \\ 
=\frac{1}{2}\left( t^{2H}-1\right) -H\ln \left( t\right) \text{ .}%
\end{array}%
\right.  \label{M45}
\end{equation}%
Eq. (\ref{M45}) proves Eq. (\ref{71}).

\subsubsection{FIM case}

For FIM, Eq. (\ref{32}) implies that 
\begin{equation}
\left. 
\begin{array}{l}
\frac{f_{t}\left( x\right) }{f_{1}\left( x\right) }=\frac{\frac{1}{2H\Gamma
\left( 1-H\right) }\cdot \left( \frac{2H^{2}}{t}\right) ^{1-H}\exp \left( -%
\frac{2H^{2}}{t}\left\vert x\right\vert ^{\frac{1}{H}}\right) \left\vert
x\right\vert ^{\frac{1}{H}-2}}{\frac{1}{2H\Gamma \left( 1-H\right) }\left(
2H^{2}\right) ^{1-H}\exp \left( -2H^{2}\left\vert x\right\vert ^{\frac{1}{H}%
}\right) \left\vert x\right\vert ^{\frac{1}{H}-2}} \\ 
\text{ } \\ 
=\frac{1}{t^{1-H}}\exp \left[ 2H^{2}\left\vert x\right\vert ^{\frac{1}{H}%
}\left( 1-\frac{1}{t}\right) \right] \text{ ,}%
\end{array}%
\right.   \label{M51}
\end{equation}%
and hence%
\begin{equation}
\ln \left[ \frac{f_{t}\left( x\right) }{f_{1}\left( x\right) }\right]
=-\left( 1-H\right) \ln \left( t\right) +2H^{2}\left( 1-\frac{1}{t}\right)
\cdot \left\vert x\right\vert ^{\frac{1}{H}}\text{ .}  \label{M52}
\end{equation}%
Setting the random motion $X\left( t\right) $ ($t\geq 0$) to be FIM, and
substituting Eq. (\ref{M52}) into Eq. (\ref{M40}) yields%
\begin{equation}
\left. 
\begin{array}{l}
\mathcal{D}\left[ I_{H}\left( t\right) ||I_{H}\left( 1\right) \right]  \\ 
\text{ } \\ 
=\int_{-\infty }^{\infty }\left[ -\left( 1-H\right) \ln \left( t\right)
+2H^{2}\left( 1-\frac{1}{t}\right) \cdot \left\vert x\right\vert ^{\frac{1}{H%
}}\right] f_{t}\left( x\right) dx \\ 
\text{ } \\ 
=-\left( 1-H\right) \ln \left( t\right) \int_{-\infty }^{\infty }f_{t}\left(
x\right) dx+2H^{2}\left( 1-\frac{1}{t}\right) \int_{-\infty }^{\infty
}\left\vert x\right\vert ^{\frac{1}{H}}f_{t}\left( x\right) dx \\ 
\text{ } \\ 
=-\left( 1-H\right) \ln \left( t\right) +2H^{2}\left( 1-\frac{1}{t}\right) 
\mathbf{E}\left[ \left\vert I_{H}\left( t\right) \right\vert ^{\frac{1}{H}}%
\right] \text{ .}%
\end{array}%
\right.   \label{M53}
\end{equation}%
The selfsimilarity of FIM implies that $I_{H}\left( t\right) =t^{H}\cdot
I_{H}\left( 1\right) $ (the equality being in law), and hence%
\begin{equation}
\left. 
\begin{array}{l}
\mathbf{E}\left[ \left\vert I_{H}\left( t\right) \right\vert ^{\frac{1}{H}}%
\right] =\mathbf{E}\left[ \left\vert t^{H}\cdot I_{H}\left( 1\right)
\right\vert ^{\frac{1}{H}}\right]  \\ 
\text{ } \\ 
=\mathbf{E}\left[ t\cdot \left\vert I_{H}\left( 1\right) \right\vert ^{\frac{%
1}{H}}\right] =t\cdot \mathbf{E}\left[ \left\vert I_{H}\left( 1\right)
\right\vert ^{\frac{1}{H}}\right] \text{ .}%
\end{array}%
\right.   \label{M54}
\end{equation}%
Using Eq. (\ref{32}), and setting $\lambda =2H^{2}$, note that the density
of the random variable $\left\vert I_{H}\left( 1\right) \right\vert $ is%
\begin{equation}
\frac{\lambda ^{1-H}}{H\Gamma \left( 1-H\right) }\exp \left( -\lambda x^{%
\frac{1}{H}}\right) x^{\frac{1}{H}-2}  \label{M55}
\end{equation}%
($0<x<\infty $), and hence (using the change-of-variables $y=\lambda x^{%
\frac{1}{H}}$):%
\begin{equation}
\left. 
\begin{array}{l}
\mathbf{E}\left[ \left\vert I_{H}\left( 1\right) \right\vert ^{\frac{1}{H}}%
\right]  \\ 
\text{ } \\ 
=\int_{0}^{\infty }x^{\frac{1}{H}}\cdot \frac{\lambda ^{1-H}}{H\Gamma \left(
1-H\right) }\exp \left( -\lambda x^{\frac{1}{H}}\right) x^{\frac{1}{H}-2}dx
\\ 
\text{ } \\ 
=\frac{\lambda ^{1-H}}{H\Gamma \left( 1-H\right) }\int_{0}^{\infty }\exp
\left( -\lambda x^{\frac{1}{H}}\right) x^{\frac{2}{H}-2}dx \\ 
\text{ } \\ 
=\frac{1}{\lambda \Gamma \left( 1-H\right) }\int_{0}^{\infty }\exp \left(
-y\right) y^{(2-H)-1}dy \\ 
\text{ } \\ 
=\frac{1}{\lambda \Gamma \left( 1-H\right) }\Gamma \left( 2-H\right)  \\ 
\text{ } \\ 
=\frac{1-H}{\lambda }=\frac{1-H}{2H^{2}}\text{ .}%
\end{array}%
\right.   \label{M56}
\end{equation}%
Combining Eq. (\ref{M53}) and Eq. (\ref{M54}) together with Eq. (\ref{M56}),
we obtain that 
\begin{equation}
\left. 
\begin{array}{l}
\mathcal{D}\left[ I_{H}\left( t\right) ||I_{H}\left( 1\right) \right]  \\ 
\text{ } \\ 
=-\left( 1-H\right) \ln \left( t\right) +2H^{2}\left( 1-\frac{1}{t}\right) 
\mathbf{E}\left[ \left\vert I_{H}\left( t\right) \right\vert ^{\frac{1}{H}}%
\right]  \\ 
\text{ } \\ 
-\left( 1-H\right) \ln \left( t\right) +2H^{2}\left( 1-\frac{1}{t}\right)
\cdot t\frac{1-H}{2H^{2}} \\ 
\text{ } \\ 
=\left( 1-H\right) \left[ t-1-\ln \left( t\right) \right] \text{ .}%
\end{array}%
\right.   \label{M57}
\end{equation}%
Eq. (\ref{M57}) proves Eq. (\ref{72}).

\subsection{The mean and the variance of the FIM increment}

Consider the FIM increment $I_{H}\left( t+\Delta \right) -I_{H}\left(
t\right) $. As FIM\ is a symmetric process, its positions have zero means,
and hence so do its increments; thus, in particular, the increment $%
I_{H}\left( t+\Delta \right) -I_{H}\left( t\right) $ has a zero mean. The
variance of the increment $I_{H}\left( t+\Delta \right) -I_{H}\left(
t\right) $ satisfies%
\begin{equation}
\left. 
\begin{array}{l}
\mathbf{Var}\left[ I_{H}\left( t+\Delta \right) -I_{H}\left( t\right) \right]
\\ 
\text{ } \\ 
=\mathbf{Var}\left[ I_{H}\left( t+\Delta \right) \right] -2\mathbf{Cov}\left[
I_{H}\left( t+\Delta \right) ,I_{H}\left( t\right) \right] +\mathbf{Var}%
\left[ I_{H}\left( t\right) \right] \text{ .}%
\end{array}%
\right.  \label{M30}
\end{equation}

As FIM is a selfsimilar process with Hurst exponent $H$, the random variable 
$I_{H}\left( t\right) $ is equal in law to the random variable $t^{H}\cdot
I_{H}\left( 1\right) $, and the random variable $I_{H}\left( t+\Delta
\right) $ is equal in law to the random variable $\left( t+\Delta \right)
^{H}\cdot I_{H}\left( 1\right) $. Hence, setting $v_{1}=\mathbf{Var}\left[
I_{H}\left( 1\right) \right] $, we have%
\begin{equation}
\mathbf{Var}\left[ I_{H}\left( t\right) \right] =v_{1}\cdot t^{2H}\text{ ,}
\label{M31}
\end{equation}%
and 
\begin{equation}
\mathbf{Var}\left[ I_{H}\left( t+\Delta \right) \right] =v_{1}\cdot \left(
t+\Delta \right) ^{2H}\text{ .}  \label{M32}
\end{equation}

The fact that the positions of FIM\ have zero means implies that%
\begin{equation}
\mathbf{Cov}\left[ I_{H}\left( t+\Delta \right) ,I_{H}\left( t\right) \right]
=\mathbf{E}\left[ I_{H}\left( t+\Delta \right) \cdot I_{H}\left( t\right) %
\right] \text{ .}  \label{M33}
\end{equation}
Denote by $\mathcal{F}_{t}$ the sigma-field generated by the trajectory of
FIM over the temporal interval $\left[ 0,t\right] $. The martingale property
of FIM implies that -- given the information $\mathcal{F}_{t}$ -- the
conditional mean of the random variable $I_{H}\left( t+\Delta \right) $ is: 
$\mathbf{E}[I_{H}\left( t+\Delta \right) | \mathcal{F}_{t}]=I_{H}\left( t\right)$. Hence, applying conditioning to the right-hand
side of Eq. (\ref{M33}), and using the fact that the random variable 
$I_{H}\left( t\right) $ has a zero mean, we have
\begin{equation}
\left. 
\begin{array}{l}
\mathbf{E}\left[ I_{H}\left( t+\Delta \right) \cdot I_{H}\left( t\right) %
\right] =\mathbf{E}\left[ \mathbf{E}\left[ I_{H}\left( t+\Delta \right)
\cdot I_{H}\left( t\right) | \mathcal{F}_{t}\right] %
\right] \\ 
\text{ } \\ 
=\mathbf{E}\left[ I_{H}\left( t\right) \cdot \mathbf{E}\left[ I_{H}\left(
t+\Delta \right) | \mathcal{F}_{t}\right] \right] =%
\mathbf{E}\left[ I_{H}\left( t\right) \cdot I_{H}\left( t\right) \right] \\ 
\text{ } \\ 
=\mathbf{E}\left[ I_{H}\left( t\right) ^{2}\right] =\mathbf{Var}\left[
I_{H}\left( t\right) \right] \text{ .}%
\end{array}%
\right.  \label{M34}
\end{equation}

Combining together Eq. (\ref{M31}), Eq. (\ref{M33}), and Eq. (\ref{M34})
yields%
\begin{equation}
\mathbf{Cov}\left[ I_{H}\left( t+\Delta \right) ,I_{H}\left( t\right) \right]
=\mathbf{Var}\left[ I_{H}\left( t\right) \right] =v_{1}\cdot t^{2H}\text{ .}
\label{M35}
\end{equation}%
In turn, substituting Eqs. (\ref{M31})-(\ref{M32}) and Eq. (\ref{M35}) into
Eq. (\ref{M30}) yields%
\begin{equation}
\left. 
\begin{array}{l}
\mathbf{Var}\left[ I_{H}\left( t+\Delta \right) -I_{H}\left( t\right) \right]
\\ 
\text{ } \\ 
=v_{1}\cdot \left( t+\Delta \right) ^{2H}-2v_{1}\cdot t^{2H}+v_{1}\cdot t^{2H} \\ 
\text{ } \\ 
=v_{1}\cdot \left[ \left( t+\Delta \right) ^{2H}-t^{2H}\right] \text{ .}%
\end{array}%
\right.  \label{M36}
\end{equation}%
Eq. (\ref{M36}) proves Eq. (\ref{52}).

In the limit $t\rightarrow \infty $ note that%
\begin{equation}
\lim_{t\rightarrow \infty }\frac{\left( t+\Delta \right) ^{2H}-t^{2H}}{%
t^{2H-1}}=\Delta \cdot \lim_{t\rightarrow \infty }\frac{\left( 1+\frac{%
\Delta }{t}\right) ^{2H}-1}{\frac{\Delta }{t}}=\Delta \cdot 2H\text{ .}
\label{M37}
\end{equation}%
Hence, in the limit $t\rightarrow \infty $, Eq. (\ref{M36}) and Eq. (\ref%
{M37}) yield the following asymptotic equivalence:%
\begin{equation}
\mathbf{Var}\left[ I_{H}\left( t+\Delta \right) -I_{H}\left( t\right) \right]
\approx \left( v_{1}\Delta 2H\right) \cdot t^{2H-1}\text{ .}  \label{M38}
\end{equation}

\subsection{The linkage between FIM and Langevin stochastic dynamics}

The stochastic dynamics of FIM are governed by the Ito SDE $\dot{I}%
_{H}\left( t\right) =\sigma \left[ I_{H}\left( t\right) \right] \dot{B}%
\left( t\right) $, with the power-law \textquotedblleft
volatility\textquotedblright\ $\sigma \left( x\right) =\left\vert
x\right\vert ^{1-\frac{1}{2H}}$. Consider a monotone increasing function $%
\varphi \left( x\right) $ that maps the real line to the real line, and
introduce the following map of FIM: the random motion $\xi _{H}\left(
t\right) =\varphi \lbrack I_{H}\left( t\right) ]$ ($t\geq 0$). Ito's formula 
\cite{KaS} and the stochastic dynamics of FIM imply that%
\begin{equation}
\dot{\xi}_{H}\left( t\right) =\frac{1}{2}\varphi ^{\prime \prime
}[I_{H}\left( t\right) ]\sigma \left[ I_{H}\left( t\right) \right]
^{2}+\left\{ \varphi ^{\prime }[I_{H}\left( t\right) ]\sigma \left[
I_{H}\left( t\right) \right] \right\} \dot{B}\left( t\right) \text{ .}
\label{M60}
\end{equation}%
Our goal is that the stochastic dynamics of random motion $\xi _{H}\left(
t\right) $ ($t\geq 0$) be governed by a Langevin SDE of the form $\dot{\xi}%
_{H}\left( t\right) =\mu \left[ \dot{\xi}_{H}\left( t\right) \right] +\dot{B}%
\left( t\right) $. To that end, the function $\varphi \left( x\right) $ must
satisfy the condition $\varphi ^{\prime }\left( x\right) \sigma \left(
x\right) =1$. As $\sigma \left( x\right) =\left\vert x\right\vert ^{1-\frac{1%
}{2H}}$, the only function $\varphi \left( x\right) $ that satisfies this
condition is 
\begin{equation}
\varphi \left( x\right) =2H\left\vert x\right\vert ^{\frac{1}{2H}}\cdot
sign\left( x\right) \text{ ,}  \label{M61}
\end{equation}%
where $sign\left( x\right) $ is the sign of $x$. Indeed, the first
derivative of the function appearing in Eq. (\ref{M61}) is%
\begin{equation}
\varphi ^{\prime }\left( x\right) =\left\vert x\right\vert ^{\frac{1}{2H}-1}=%
\frac{1}{\sigma \left( x\right) }\text{ .}  \label{M62}
\end{equation}%
Also, the second derivative of the function appearing in Eq. (\ref{M61}) is 
\begin{equation}
\varphi ^{\prime \prime }\left( x\right) =\left( \frac{1}{2H}-1\right)
\left\vert x\right\vert ^{\frac{1}{2H}-2}\cdot sign\left( x\right) \text{ .}
\label{M63}
\end{equation}%
Eqs. (\ref{M61})-(\ref{M63}) imply that%
\begin{equation}
\left. 
\begin{array}{l}
\frac{1}{2}\varphi ^{\prime \prime }\left( x\right) \sigma \left( x\right)
^{2} \\ 
\text{ } \\ 
=\frac{1}{2}\left( \frac{1}{2H}-1\right) \left\vert x\right\vert ^{\frac{1}{%
2H}-2}\cdot sign\left( x\right) \cdot \left\vert x\right\vert ^{2-\frac{1}{H}%
} \\ 
\text{ } \\ 
=\frac{1}{2}\left( \frac{1}{2H}-1\right) \left\vert x\right\vert ^{-\frac{1}{%
2H}}\cdot sign\left( x\right)  \\ 
\text{ } \\ 
=\frac{1}{2}\left( 1-2H\right) \frac{1}{2H\left\vert x\right\vert ^{\frac{1}{%
2H}}\cdot sign\left( x\right) } \\ 
\text{ } \\ 
=\left( \frac{1}{2}-H\right) \frac{1}{\varphi \left( x\right) }\text{ .}%
\end{array}%
\right.   \label{M64}
\end{equation}%
In turn, substituting Eq. (\ref{M62}) and Eq. (\ref{M64}) into Eq. (\ref{M60}%
) implies that 
\begin{equation}
\left. 
\begin{array}{l}
\dot{\xi}_{H}\left( t\right) =\left( \frac{1}{2}-H\right) \frac{1}{\varphi
\lbrack I_{H}\left( t\right) ]}+\dot{B}\left( t\right)  \\ 
\text{ } \\ 
=\left( \frac{1}{2}-H\right) \frac{1}{\xi _{H}\left( t\right) }+\dot{B}%
\left( t\right) \text{ .}%
\end{array}%
\right.   \label{M65}
\end{equation}%
Eq. (\ref{M65}) proves Eq. (\ref{80}).

\bigskip\


\begin{thebibliography}{999}
\bibitem{SSDIF} Eliazar, Iddo Isaac. "Selfsimilar diffusions." Journal of
Physics A: Mathematical and Theoretical (2021).

\bibitem{EM} Embrechts, Paul, and Makoto Maejima. Selfsimilar processes.
Princeton University Press, 2009.

\bibitem{Fri} Friedman,Avner. Stochastic differential equations and
applications. Dover publications, 2006.

\bibitem{Oks} Oksendal, Bernt. Stochastic differential equations: an
introduction with applications. Springer Science \& Business Media, 2013.

\bibitem{Arn} Arnold, Ludwig. Stochastic differential equations: Theory and
applications. Dover publications, 2014.\ 

\bibitem{PSW} Pekalski, Andrzej, and Katarzyna Sznajd-Weron (Eds.).
Anomalous Diffusion: From Basics to Applications. Springer, 1999.

\bibitem{KSb} Klafter, Joseph, and Igor M. Sokolov. First steps in random
walks: from tools to applications. Oxford University Press, 2011.

\bibitem{MCB} Mendez, Vicenc, Daniel Campos, and Frederic Bartumeus.
Stochastic foundations in movement ecology: Anomalous Diffusion, Front
Propagation and Random Searches. Springer-Verlag, 2013.

\bibitem{EKL} Evangelista, Luiz Roberto, and Ervin Kaminski Lenzi.
Fractional diffusion equations and anomalous diffusion. Cambridge University
Press, 2018.

\bibitem{LMO} Lindenberg, Katja, Ralf Metzler, and Gleb Oshanin, eds.
Chemical kinetics: beyond the textbook. World scientific, 2019.

\bibitem{DHWX} Deng, Weihua, Ru Hou, Wanli Wang, and Pengbo Xu. Modeling
anomalous diffusion: from statistics to mathematics. World Scientific, 2020.

\bibitem{Shl} Shlesinger, Michael F. An Unbounded Experience in Random Walks
with Applications. World Scientific, 2021.

\bibitem{Gar} Gardiner, Crispin W. Handbook of stochastic methods. Springer,
1985.

\bibitem{Kam} Van Kampen, Nicolaas Godfried. Stochastic processes in physics
and chemistry. North Holland, 2007.

\bibitem{KaS} Karatzas, Ioannis, and Steven E. Shreve. Brownian motion and
stochastic calculus. Springer, 1991.

\bibitem{CheG} Chechkin, A. V., and V. Yu Gonchar. "Fractional Brownian
motion approximation based on fractional integration of a white noise."
Chaos, Solitons \& Fractals 12, no. 2 (2001): 391-398.

\bibitem{MKW} Mercik, Sz. and K. Weron. "Fractional Brownian Motion as a
Model of the Self-Similar Ion Channel Kinetics." Acta Physica Polonica.
Series B 32, no. 5 (2001): 1621-1630.

\bibitem{MWBK} Magdziarz, Marcin, Aleksander Weron, Krzysztof Burnecki, and
Joseph Klafter. "Fractional Brownian motion versus the continuous-time
random walk: A simple test for subdiffusive dynamics." Physical review
letters 103, no. 18 (2009): 180602.

\bibitem{JeoM} Jeon, Jae-Hyung, and Ralf Metzler. "Fractional Brownian
motion and motion governed by the fractional Langevin equation in confined
geometries." Physical Review E 81, no. 2 (2010): 021103.

\bibitem{JCM1} Jeon, Jae-Hyung, A. V. Chechkin, and Ralf Metzler. "First passage
behaviour of fractional Brownian motion in two-dimensional wedge domains."
EPL (Europhysics Letters) 94, no. 2 (2011): 20008.

\bibitem{BKJB} Burnecki, Krzysztof, Eldad Kepten, Joanna Janczura, Irena
Bronshtein, Yuval Garini, and Aleksander Weron. "Universal algorithm for
identification of fractional Brownian motion. A case of telomere
subdiffusion." Biophysical journal 103, no. 9 (2012): 1839-1847.

\bibitem{JCM2} Jeon, Jae-Hyung, Aleksei V. Chechkin, and Ralf Metzler.
"First passage behavior of multi-dimensional fractional Brownian motion and
application to reaction phenomena." In First-Passage Phenomena and Their
Applications, pp. 175-202. 2014.

\bibitem{SBW} Sikora, Grzegorz, Krzysztof Burnecki, and Agnieszka Wy\l oma%
\'{n}ska. "Mean-squared-displacement statistical test for fractional
Brownian motion." Physical Review E 95, no. 3 (2017): 032110.

\bibitem{GPV} Guggenberger, Tobias, Gianni Pagnini, Thomas Vojta, and Ralf
Metzler. "Fractional Brownian motion in a finite interval: correlations
effect depletion or accretion zones of particles near boundaries." New
Journal of Physics 21, no. 2 (2019): 022002.

\bibitem{WCCTS} Wang, Wei, Andrey G. Cherstvy, Aleksei V. Chechkin,
Samudrajit Thapa, Flavio Seno, Xianbin Liu, and Ralf Metzler. "Fractional
Brownian motion with random diffusivity: emerging residual nonergodicity
below the correlation time." Journal of Physics A: Mathematical and
Theoretical 53, no. 47 (2020): 474001.

\bibitem{GCM} Guggenberger, Tobias, Aleksei Chechkin, and Ralf Metzler.
"Fractional Brownian motion in superharmonic potentials and non-Boltzmann
stationary distributions." Journal of Physics A: Mathematical and
Theoretical 54, no. 29 (2021): 29LT01.

\bibitem{WCKM} Wang, Wei, Andrey Cherstvy, Holger Kantz, Ralf Metzler, and
Igor Sokolov. "Time-averaging and emerging nonergodicity upon resetting of
fractional Brownian motion and heterogeneous diffusion processes." bioRxiv
(2021).

\bibitem{CWMS} Cherstvy, Andrey, Wei Wang, Ralf Metzler, and Igor Sokolov.
"Inertia triggers nonergodicity of fractional Brownian motion." bioRxiv
(2021).

\bibitem{CCM} Cherstvy, Andrey G., Aleksei V. Chechkin, and Ralf Metzler. "Anomalous diffusion and ergodicity breaking in heterogeneous diffusion processes." New Journal of Physics 15, no. 8 (2013): 083039.

\bibitem{LB} Leibovich, N., and E. Barkai. "Infinite ergodic theory for heterogeneous diffusion processes." Physical Review E 99, no. 4 (2019): 042138.

\bibitem{LLG} Lenzi, M. K., E. K. Lenzi, L. M. S. Guilherme, L. R. Evangelista, and H. V. Ribeiro. "Transient anomalous diffusion in heterogeneous media with stochastic resetting." Physica A: Statistical Mechanics and its Applications (2021): 126560.

\bibitem{Vpas} Van der Pas, Peter W. "The discovery of the Brownian motion."
Scientiarum Historia: Tijdschriftvoor de Geschiedenis van de Wetenschappenen
de Geneeskunde 13, no. 1 (1971): 27-35.

\bibitem{Bro} Brown, Robert. "XXVII. A brief account of microscopical
observations made in the months of June, July and August 1827, on the
particles contained in the pollen of plants; and on the general existence of
active molecules in organic and inorganic bodies." The philosophical
magazine 4, no. 21 (1828): 161-173.

\bibitem{Per} Perrin, Jean. "Mouvement Brownien et realite moleculaire."
Annales de Chimie et de Physique 18 (1909): 5-104.

\bibitem{SL1} Scher, H., and M. Lax. "Continuous time random walk model of
hopping transport: application to impurity conduction." Journal of
Non-Crystalline Solids 8 (1972): 497-504.

\bibitem{SL2} Scher, Harvey, and Melvin Lax. "Stochastic transport in a
disordered solid. I. Theory." Physical Review B 7, no. 10 (1973): 4491.

\bibitem{SL3} Scher, Harvey, and Melvin Lax. "Stochastic transport in a
disordered solid. II. Impurity conduction." Physical Review B 7, no. 10
(1973): 4502.

\bibitem{Shl1} Shlesinger, Michael F. "Asymptotic solutions of
continuous-time random walks." Journal of Statistical Physics 10, no. 5
(1974): 421-434.

\bibitem{SMon} Scher, Harvey, and Elliott W. Montroll. "Anomalous
transit-time dispersion in amorphous solids." Physical Review B 12, no. 6
(1975): 2455.

\bibitem{GAA} Gefen, Yuval, Amnon Aharony, and Shlomo Alexander. "Anomalous
diffusion on percolating clusters." Physical Review Letters 50, no. 1
(1983): 77.

\bibitem{KBS} Klafter, Joseph, Alexander Blumen, and Michael F. Shlesinger.
"Stochastic pathway to anomalous diffusion." Physical Review A 35, no. 7
(1987): 3081.

\bibitem{BG} Bouchaud, Jean-Philippe, and Antoine Georges. "Anomalous
diffusion in disordered media: statistical mechanisms, models and physical
applications." Physics reports 195, no. 4-5 (1990): 127-293.

\bibitem{Zan} Zanette, Damian H., and Pablo A. Alemany. "Thermodynamics of
anomalous diffusion." Physical review letters 75, no. 3 (1995): 366.

\bibitem{MK} Metzler, Ralf, and Joseph Klafter. "The random walk's guide to
anomalous diffusion: a fractional dynamics approach." Physics reports 339,
no. 1 (2000): 1-77.

\bibitem{SK} Sokolov, Igor M., and Joseph Klafter. "From diffusion to
anomalous diffusion: a century after Einstein's Brownian motion." Chaos: An
Interdisciplinary Journal of Nonlinear Science 15, no. 2 (2005): 026103.

\bibitem{KS} Klafter, Joseph, and Igor M. Sokolov. "Anomalous diffusion
spreads its wings." Physics world 18, no. 8 (2005): 29.

\bibitem{Sok} Sokolov, Igor M. "Models of anomalous diffusion in crowded
environments." Soft Matter 8, no. 35 (2012): 9043-9052.

\bibitem{MJCB} Metzler, Ralf, Jae-Hyung Jeon, Andrey G. Cherstvy, and Eli
Barkai. "Anomalous diffusion models and their properties: non-stationarity,
non-ergodicity, and ageing at the centenary of single particle tracking."
Physical Chemistry Chemical Physics 16, no. 44 (2014): 24128-24164.

\bibitem{San} dos Santos, Maike AF. "Analytic approaches of the anomalous
diffusion: A review." Chaos, Solitons \& Fractals 124 (2019): 86-96.

\bibitem{Metz} Metzler, Ralf. "Brownian motion and beyond: first-passage,
power spectrum, non-Gaussianity, and anomalous diffusion." Journal of
Statistical Mechanics: Theory and Experiment 2019, no. 11 (2019): 114003.

\bibitem{SLLS} Sancho, Jose Maria, A. M. Lacasta, Katja Lindenberg, Igor M.
Sokolov, and A. H. Romero. "Diffusion on a solid surface: Anomalous is
normal." Physical review letters 92, no. 25 (2004): 250601.

\bibitem{AnoUbi} Eliazar, Iddo, and Joseph Klafter. "Anomalous is
ubiquitous." Annals of Physics 326, no. 9 (2011): 2517-2531.

\bibitem{VSSP} Vitali, S., Sposini, V., Sliusarenko, O., Paradisi, P.,
Castellani, G., \& Pagnini, G. (2018). Langevin equation in complex media
and anomalous diffusion. Journal of The Royal Society Interface, 15(145),
20180282.

\bibitem{OFL} Oliveira, Fernando A., Rogelma Ferreira, Luciano C. Lapas, and
Mendeli H. Vainstein. "Anomalous diffusion: A basic mechanism for the
evolution of inhomogeneous systems." Frontiers in Physics 7 (2019): 18.

\bibitem{SKla} Sato, Yuzuru, and Rainer Klages. "Anomalous diffusion in
random dynamical systems." Physical review letters 122, no. 17 (2019):
174101.

\bibitem{BSFV} Bo, Stefano, Falko Schmidt, Ralf Eichhorn, and Giovanni
Volpe. "Measurement of anomalous diffusion using recurrent neural networks."
Physical Review E 100, no. 1 (2019): 010102.

\bibitem{WCLM} Wang, Wei, Andrey G. Cherstvy, Xianbin Liu, and Ralf Metzler.
"Anomalous diffusion and nonergodicity for heterogeneous diffusion processes
with fractional Gaussian noise." Physical Review E 102, no. 1 (2020): 012146.

\bibitem{FTPB} Feng, Libo, Ian Turner, Patrick Perre, and Kevin Burrage. "An
investigation of nonlinear time-fractional anomalous diffusion models for
simulating transport processes in heterogeneous binary media."
Communications in Nonlinear Science and Numerical Simulation 92 (2021):
105454.

\bibitem{CSRM} Cherstvy, Andrey G., Hadiseh Safdari, and Ralf Metzler.
"Anomalous diffusion, nonergodicity, and ageing for exponentially and
logarithmically time-dependent diffusivity: striking differences for massive
versus massless particles." Journal of Physics D: Applied Physics 54, no. 19
(2021): 195401.

\bibitem{Kon} Kononovicius, A. "Noisy voter model for the anomalous
diffusion of parliamentary presence." Journal of Statistical Mechanics:
Theory and Experiment 2020, no. 6 (2020): 063405.

\bibitem{KazK} Kazakevi\v{c}ius, Rytis, and Aleksejus Kononovicius.
"Anomalous diffusion in nonlinear transformations of the noisy voter model."
Physical Review E 103, no. 3 (2021): 032154.

\bibitem{MasM} Masoliver, Jaume, and Miquel Montero. "Anomalous diffusion
under stochastic resettings: A general approach." Physical Review E 100, no.
4 (2019): 042103.

\bibitem{FauS} Antonio Faustino dos Santos, Maike. "Comb model with
non-static stochastic resetting and anomalous diffusion." Fractal and
Fractional 4, no. 2 (2020): 28.

\bibitem{DzS} Dzhanoev, A. R., and I. M. Sokolov. "The effect of the
junction model on the anomalous diffusion in the 3D comb structure." Chaos,
Solitons \& Fractals 106 (2018): 330-336.

\bibitem{LSRJ} Lenzi, E. K., T. Sandev, H. V. Ribeiro, P. Jovanovski, A. Iomin, and L. Kocarev. "Anomalous diffusion and random search in xyz-comb:
exact results." Journal of Statistical Mechanics: Theory and Experiment
2020, no. 5 (2020): 053203.

\bibitem{Faus} Antonio Faustino dos Santos, Maike. "Comb model with
non-static stochastic resetting and anomalous diffusion." Fractal and
Fractional 4, no. 2 (2020): 28.

\bibitem{Iom} Iomin, A. "Anomalous diffusion in umbrella comb." Chaos,
Solitons \& Fractals 142 (2021): 110488.

\bibitem{WSSC} Wang, Wei, Flavio Seno, Igor M. Sokolov, Aleksei V. Chechkin,
and Ralf Metzler. "Unexpected crossovers in correlated random-diffusivity
processes." New Journal of Physics 22, no. 8 (2020): 083041.

\bibitem{WCCT} @ Wang, Wei, Andrey G. Cherstvy, Aleksei V. Chechkin,
Samudrajit Thapa, Flavio Seno, Xianbin Liu, and Ralf Metzler. "Fractional
Brownian motion with random diffusivity: emerging residual nonergodicity
below the correlation time." Journal of Physics A: Mathematical and
Theoretical 53, no. 47 (2020): 474001.

\bibitem{GSMO} Grebenkov, Denis S., Vittoria Sposini, Ralf Metzler, Gleb
Oshanin, and Flavio Seno. "Exact first-passage time distributions for three
random diffusivity models." Journal of Physics A: Mathematical and
Theoretical 54, no. 4 (2021): 04LT01.

\bibitem{SanM} dos Santos, Maike AF, and Luiz Menon Junior. "Random
diffusivity models for scaled Brownian motion." Chaos, Solitons \& Fractals
144 (2021): 110634.

\bibitem{Bac1} Bachelier, Louis. "Theorie de la speculation." In Annales
scientifiques de lecole normale superieure, vol. 17, pp. 21-86. 1900.

\bibitem{Bac2} Bachelier, Louis. Louis Bachelier's theory of speculation:
the origins of modern finance. Princeton University Press, 2011.

\bibitem{Ein} Einstein, Albert. "Uber die von der molekularkinetischen
theorie der warmegeforderte bewegung von in ruhenden
flussigkeitensuspendierten teilchen." Annalen der physik 4 (1905).

\bibitem{Smo} Von Smoluchowski, Marian. "Zurkinetischentheorie der
brownschen molekularbewegung und der suspensionen." Annalen der physik 326,
no. 14 (1906): 756-780.

\bibitem{Wie} Wiener, Norbert. "Differential-Space." Journal of Mathematics
and Physics 2, no. 1-4 (1923): 131-174.

\bibitem{Kol} Kolmogorov, Andrei Nikolaevitch. "The Wiener spiral and some
other interesting curves in Hilbert space." In Dokl. Akad. Nauk SSSR, vol.
26, no. 2 (1940): 115-118.

\bibitem{Yag} Yaglom, A. M. "Correlation theory of processes with
stationary increments of order n" Amer. Math. Soc. Transl. Series, American
Math. Soc. Providence RI 8 (1958): 87-102.

\bibitem{MV} Mandelbrot, Benoit B., and John W. Van Ness. "Fractional
Brownian motions, fractional noises and applications." SIAM review 10, no. 4
(1968): 422-437.

\bibitem{Mish} Yuliya Mishura. Stochastic calculus for fractional Brownian
motion and related processes. Springer Science \& Business Media, 2008.

\bibitem{BHOZ} Biagini, Francesca, Yaozhong Hu, Bernt Oksendal, and Tusheng
Zhang. Stochastic calculus for fractional Brownian motion and applications.
Springer Science \& Business Media, 2008.

\bibitem{Nour} Nourdin, Ivan. Selected aspects of fractional Brownian
motion. Springer, 2012.

\bibitem{BMRS} Banna, Oksana, Yuliya Mishura, Kostiantyn Ralchenko, and
Sergiy Shklyar. Fractional Brownian Motion: Approximations and Projections.
John Wiley \& Sons, 2019.

\bibitem{FracMo} Eliazar, Iddo I., and Michael F. Shlesinger. "Fractional
motions." Physics Reports 527, no. 2 (2013): 101-129.

\bibitem{Yin} Yin, Z-M. "New methods for simulation of fractional Brownian
motion." Journal of computational physics 127, no. 1 (1996): 66-72.

\bibitem{Coe} Jean-Francois, Coeurjolly. "Simulation and identification of
the fractional Brownian motion: a bibliographical and comparative study."
Journal of statistical software 5 (2000): 1-53.

\bibitem{Ste} Stein, Michael L. "Fast and exact simulation of fractional
Brownian surfaces." Journal of Computational and Graphical Statistics 11,
no. 3 (2002): 587-599.

\bibitem{DMen} Dieker, Antonius Bernardus, and Michael Mandjes. "On spectral
simulation of fractional Brownian motion." Probability in the Engineering
and Informational Sciences 17, no. 3 (2003): 417-434.

\bibitem{Diek} Dieker, Ton. "Simulation of fractional Brownian motion." PhD
diss., Masters Thesis, Department of Mathematical Sciences, University of
Twente, The Netherlands, 2004.

\bibitem{Pip} Pipiras, Vladas. "Wavelet-based simulation of fractional
Brownian motion revisited." Applied and Computational Harmonic Analysis 19,
no. 1 (2005): 49-60.

\bibitem{Pas} Pashko, Anatolii. "Accuracy of simulation of fractional
Brownian motion." Information Technologies, Management and Society (2018).

\bibitem{PSO} Pashko, Anatolii, Olga Sinyavska, and Tetiana Oleshko.
"Simulation of Fractional Brownian Motion and Estimation of Hurst
Parameter." In 2020 IEEE 15th International Conference on Advanced Trends in
Radioelectronics, Telecommunications and Computer Engineering (TCSET), pp.
632-637. IEEE, 2020.

\bibitem{CDH} Chen, Yi, Jing Dong, and Hao Ni. "$\epsilon $-Strong
Simulation of Fractional Brownian Motion and Related Stochastic Differential
Equations." Mathematics of Operations Research (2021).

\bibitem{SMB} Shahnazi-Pour, A., B. Parsa Moghaddam, and Afshin Babaei.
"Numerical simulation of the Hurst index of solutions of fractional
stochastic dynamical systems driven by fractional Brownian motion." Journal
of Computational and Applied Mathematics 386 (2021): 113210.

\bibitem{Lan} Langevin, Paul. "Sur la theorie du mouvement Brownien." Compt.
Rendus 146 (1908): 530-533.

\bibitem{CK} Coffey, William, and Yuri P. Kalmykov. The Langevin equation:
with applications to stochastic problems in physics, chemistry and
electrical engineering. World Scientific, 2012.

\bibitem{Pav} Pavliotis, Grigorios A. Stochastic processes and applications:
diffusion processes, the Fokker-Planck and Langevin equations. Vol. 60.
Springer, 2014.

\bibitem{Ito} Ito, Kiyosi. On stochastic differential equations. American
Mathematical Soc., 1951.

\bibitem{IK} Ito, Kiyosi, and P. Henry Jr. Diffusion processes and their
sample paths. Springer Science \& Business Media, 2012.

\bibitem{Hul} Hull, John C. Options, futures, and other derivatives.
Pearson, 2017.

\bibitem{CSM} Chechkin, Aleksei V., Flavio Seno, Ralf Metzler, and Igor M.
Sokolov. "Brownian yet non-Gaussian diffusion: from superstatistics to
subordination of diffusing diffusivities." Physical Review X 7, no. 2
(2017): 021002.

\bibitem{TCh} Tyagi, Neha, and Binny J. Cherayil. "Non-Gaussian Brownian
diffusion in dynamically disordered thermal environments." The Journal of
Physical Chemistry B 121, no. 29 (2017): 7204-7209.

\bibitem{Met1} Metzler, Ralf. "Gaussianity fair: the riddle of anomalous yet
non-Gaussian diffusion." Biophysical journal 112, no. 3 (2017): 413.

\bibitem{SMM} Slezak, Jakub, Ralf Metzler, and Marcin Magdziarz.
"Superstatistical generalised Langevin equation: non-Gaussian viscoelastic
anomalous diffusion." New Journal of Physics 20, no. 2 (2018): 023026.

\bibitem{SCS} Sposini, Vittoria, Aleksei V. Chechkin, Flavio Seno, Gianni
Pagnini, and Ralf Metzler. "Random diffusivity from stochastic equations:
comparison of two models for Brownian yet non-Gaussian diffusion." New
Journal of Physics 20, no. 4 (2018): 043044.

\bibitem{LG1} Lanoiselee, Yann, and Denis S. Grebenkov. "A model of
non-Gaussian diffusion in heterogeneous media." Journal of Physics A:
Mathematical and Theoretical 51, no. 14 (2018): 145602.

\bibitem{LG2} Lanoiselee, Yann, and Denis S. Grebenkov. "Non-Gaussian
diffusion of mixed origins." Journal of Physics A: Mathematical and
Theoretical 52, no. 30 (2019): 304001.

\bibitem{SBM} Slezak, Jakub, Krzysztof Burnecki, and Ralf Metzler. "Random
coefficient autoregressive processes describe Brownian yet non-Gaussian
diffusion in heterogeneous systems." New Journal of Physics 21, no. 7
(2019): 073056.

\bibitem{Met2} Metzler, Ralf. "Superstatistics and non-Gaussian diffusion."
The European Physical Journal Special Topics 229, no. 5 (2020): 711-728.

\bibitem{StaEv} Eliazar, Iddo I., and Igor M. Sokolov. "Measuring
statistical evenness: A panoramic overview." Physica A: Statistical
Mechanics and its Applications 391, no. 4 (2012): 1323-1353.

\bibitem{HarIn} Eliazar, Iddo. "Harnessing inequality." Physics Reports 649
(2016): 1-29.

\bibitem{TouIn} Eliazar, Iddo. "A tour of inequality." Annals of Physics 389
(2018): 306-332.

\bibitem{Cou} P.B. Coulter, Measuring inequality: A methodological handbook
(Westview Press, Boulder, 1989).

\bibitem{BL} G. Betti and A. Lemmi (Eds.), Advances on income inequality and
concentration measures (Routledge, New-York, 2008).

\bibitem{HN} L. Hao and D.Q. Naiman, Assessing Inequality (Sage, Los
Angeles, 2010).

\bibitem{Cow} F. Cowell, Measuring inequality (Oxford University Press,
Oxford, 2011).

\bibitem{RenEq} Eliazar, Iddo. "Investigating equality: The Renyi spectrum."
Physica A: Statistical Mechanics and its Applications 481 (2017): 90-118.

\bibitem{KL} Kullback, Solomon, and Richard A. Leibler. "On information and
sufficiency." The annals of mathematical statistics 22, no. 1 (1951): 79-86.

\bibitem{Kul} Kullback, Solomon. Information theory and statistics. Dover
Publications (new edition) 1997.

\bibitem{KB} Kessler, David A., and Eli Barkai. "Infinite covariant density
for diffusion in logarithmic potentials and optical lattices." Physical
review letters 105, no. 12 (2010): 120602.

\bibitem{DLB} Dechant, A., E. Lutz, E. Barkai, and D. A. Kessler. "Solution
of the Fokker-Planck equation with a logarithmic potential." Journal of
Statistical Physics 145, no. 6 (2011): 1524-1545.

\bibitem{HMS1} Hirschberg, Ori, David Mukamel, and Gunter M. Schutz.
"Approach to equilibrium of diffusion in a logarithmic potential." Physical
Review E 84, no. 4 (2011): 041111.

\bibitem{HMS2} Hirschberg, Ori, David Mukamel, and Gunter M. Schutz.
"Diffusion in a logarithmic potential: scaling and selection in the approach
to equilibrium." Journal of Statistical Mechanics: Theory and Experiment
2012, no. 02 (2012): P02001.

\bibitem{BAK} Barkai, E., E. Aghion, and D. A. Kessler. "From the area under
the Bessel excursion to anomalous diffusion of cold atoms." Physical Review
X 4, no. 2 (2014): 021036.

\bibitem{RBH} Ryabov, Artem, Ekaterina Berestneva, and Viktor Holubec.
"Brownian motion in time-dependent logarithmic potential: exact results for
dynamics and first-passage properties." Journal of chemical physics 143, no.
11 (2015): 114117.

\bibitem{SVB} De Santo, Serena, Pablo Villegas, Raffael la Burioni, and
Miguel A. Munoz. "Simple unified view of branching process statistics:
Random walks in balanced logarithmic potentials." Physical Review E 95, no.
3 (2017): 032115.

\bibitem{RR} Ray, Somrita, and Shlomi Reuveni. "Diffusion with resetting in
a logarithmic potential." Journal of chemical physics 152 (2020): 234110.

\bibitem{OPR} Onofri, Manuele, Gaia Pozzoli, Mattia Radice, and Roberto
Artuso. "Exploring the Gillis model: a discrete approach to diffusion in
logarithmic potentials." Journal of Statistical Mechanics: Theory and
Experiment 2020, no. 11 (2020): 113201.

\bibitem{PM} Paraguassu, Pedro V., and Welles AM Morgado. "The heat
distribution in a logarithm potential." Journal of Statistical Mechanics:
Theory and Experiment 2021, no. 2 (2021): 023205.

\bibitem{Pit} Pitman, James W. "One-dimensional Brownian motion and the
three-dimensional Bessel process." Advances in Applied Probability 7, no. 3
(1975): 511-526.

\bibitem{DeL} De Long, David M. "Crossing probabilities for a square root
boundary by a Bessel process." Communications in Statistics-Theory and
Methods 10, no. 21 (1981): 2197-2213.

\bibitem{Imh} Imhof, J-P. "Density factorizations for Brownian motion,
meander and the three-dimensional Bessel process, and applications." Journal
of Applied Probability 21, no. 3 (1984): 500-510.

\bibitem{Che} Cherny, A. S. "On the strong and weak solutions of stochastic
differential equations governing Bessel processes." Stochastics: An
International Journal of Probability and Stochastic Processes 70, no. 3-4
(2000): 213-219.

\bibitem{Kat} Katori, Makoto. Bessel processes, Schramm-Loewner evolution,
and the Dyson model. Heidelberg: Springer, 2016.

\bibitem{BGT} N. H. Bingham, C. M. Goldie, and J.L. Teugels, Regular
variation (Cambridge University Press, Cambridge, 1987).
\end{thebibliography}
\end{document}